\documentclass[leqno,11pt]{article}

\usepackage{amssymb,hyperref,amsthm,hyperref}
\usepackage{euscript}
\usepackage[dvips]{graphicx}
\usepackage{flafter}
\usepackage{pstricks}
\usepackage[latin1]{inputenc}
\usepackage{amsmath}
\usepackage{amsfonts}
\usepackage{pstricks-add}
\usepackage{dsfont}
\usepackage{bm}
\usepackage{tikz}

\hypersetup{
    linktoc=page,
    linkcolor=red,          % color of internal links
    citecolor=blue,        % color of links to bibliography
    filecolor=blue,      % color of file links
    urlcolor=cyan,
    colorlinks=true           % color of external links
}

\usepackage[refpage]{nomencl}
\makenomenclature

\usepackage{makeidx}
\makeindex

\usepackage{mathrsfs} %%% Nice caligraphic fonts
\evensidemargin\oddsidemargin
\begin{document}

%%
%%%%%%%%%%%%%%%%%%%%%%%%%%% Equation numberings
\newcommand{\eqnsection}{
\renewcommand{\theequation}{\thesection.\arabic{equation}}
   \makeatletter
   \csname  @addtoreset\endcsname{equation}{section}
   \makeatother}
\eqnsection
%%%%%%%%%%%%%%%%%%%%%%%%%%%

%%%%%%%%%%%%%% Bbb characters
%%%%%%%%%%%%%% Real numbers
\def\r{{\mathbb R}}
%%%%%%%%%%%%%% Expectation
\def\e{{\mathbb E}}
%%%%%%%%%%%%%% Probability
\def\p{{\mathbb P}}
%%%%%%%%%%%%%% Law of environment
\def\P{{\mathbb P}}
\def\E{{\mathbb E}}
\def\Q{{\bf Q}}
%%%%%%%%%%%%%% Integers
\def\rc{]\!]}
\def\lc{[\![}
\def\z{{\mathbb Z}}
%%%%%%%%%%%%%% Natural numbers
\def\N{{\mathbb N}}
%%%%%%%%%%%%%% Tree
\def\T{{\mathbb T}}
%%%%%%%%%%%%%% Galton-Watson tree
\def\G{{\mathbb G}}

\def\C{{\mathbb C}}
%%%%%%%%%%%%%% "Low"
\def\L{{\mathbb L}}
%%%%%%%%%%%%%% "indicatrice"
\def\1{{\mathds{1}}}
%%%%%%%%%%%%%%
\def\deg{\chi}

\def\d{\mathtt{d}}
%%%%%%%%%%%%%%
\def\ttheta{{\bm \theta}}
%%%%%%%%%%%%%%
%%%%%%%%%%%%%% "vecteur t"
\def\t{{\bf{t}}}
\def\a{{\bf{a}}}
%%%%%%%%%%%%%%
\def\deg{\chi}
\def\B{\mathfrak{B}}

\def\M{{\mathbb{M}}}
%%%%%%%%%%%%%%
%%%%%%%%%%%%%%%% Special symbols
%%%%%%%%%%%%%% Exponential
\def\ee{e}
%%%%%%%%%%%%%% Differentiation
\def\d{\, \mathrm{d}}
%%%%%%%%%%%%%% Survival
\def\S{\mathscr{S}}
%%%%%%%%%%%%%% Binary search
\def\bs{{\tt bs}}
%%%%%%%%%%%%%%
\def\bbeta{{\bm \beta}}
%%%%%%%%%%%%%%
\def\ttheta{{\bm \theta}}
%%%%%%%%%%%%%%
%%%%%%%%%%%%%%%%%%%%%%%%%%%%%%%%%%%%%%%%%%%%
%%%%%%%%%%%%%%%%%%%%%%%%%%%%%%%%%%%%%%%%%%%%
%%%%%%%%%%%%%%%%%%%%%%%%%%%%%%%%%%%%%%%%%%%%

\newtheorem{theorem}{Theorem}[section]

\newtheorem{definition}[theorem]{Definition}
\newtheorem{Lemma}[theorem]{Lemma}
\newtheorem{Proposition}[theorem]{Proposition}
\newtheorem{Remark}[theorem]{Remark}
\newtheorem{corollary}[theorem]{Corollary}
\newtheorem{conjecture}[theorem]{Conjecture}
%%%%%%%%%%%%%%%%%%%%%%%%%%%%%%%%%%%%%%%%%%%%
%%%%%%%%%%%%%%%%%%%%%%%%%%%%%%%%%%%%%%%%%%%%
%%%%%%%%%%%%%%%%%%%%%%%%%%%%%%%%%%%%%%%%%%%%

%%%%%%%%%%%%%% Beginning of the text

\title{Continuity estimates for the complex cascade model on the phase boundary}
\date{\vspace{-5ex}}
 
 \author{ Thomas Madaule \footnote{Universit\'e Paul Sabatier, IMT, Toulouse, France} , R\'emi Rhodes \footnote{Universit{\'e} Paris-Est Marne la Vall\'ee, LAMA, Champs sur Marne, France.},
 Vincent Vargas \footnote{ENS Ulm, DMA, 45 rue d'Ulm,  75005 Paris, France.} }

 \maketitle

\begin{abstract}
We consider the complex branching random walk on a dyadic tree with Gaussian weights on the boundary between the diffuse phase and the glassy phase. We study the branching random walk in the space of continuous functions and establish convergence in this space. The main difficulty here is that the expected modulus of continuity of the limit is too weak in order to show tightness in the space of continuous functions by means of  standard tools from the theory of stochastic processes. 
\end{abstract}
 
\noindent{\bf Key words or phrases:}  branching random walk, complex weights,  boundary case, diffusivity.


{\leftskip=2truecm \rightskip=2truecm \baselineskip=15pt \small

} %%%%%% End of narrower

\tableofcontents

%%%%%%%%%%%%%%%%%%%%%%%%%%%%%%%%%%%%%%%%%%%%%%%%%%%%%%%%%%%%%%%%%%%%%%%%%%%%%%%%
\section{Introduction}
%%%%%%%%%%%%%%%%%%%%%%%%%%%%%%%%%%%%%%%%%%%%%%%%%%%%%%%%%%%%%%%%%%%%%%%%%%%%%%%%

\subsection{The model and main results}
 
We consider a discrete-time complex valued branching random walk. The system starts with an initial particle, called the root, at time $n=0$. At time $n=1$, the particle dies and gives birth to $2$ particles, which form the particles at generation $1$. At time $n=2$, each of these particles dies and gives birth to $2$ new particles, and so on... For all $n \geq 0$, we denote $\mathbb{T}_{\leq n }=\{0,1\}^{\{1,\dots,n\}} $ the genealogical tree associated to the $n$-th generation: its elements   have length $n$ which we denote by $|u|=n$. Then, we set $\mathbb{T}= \bigcup_n \mathbb{T}_{\leq n } $.

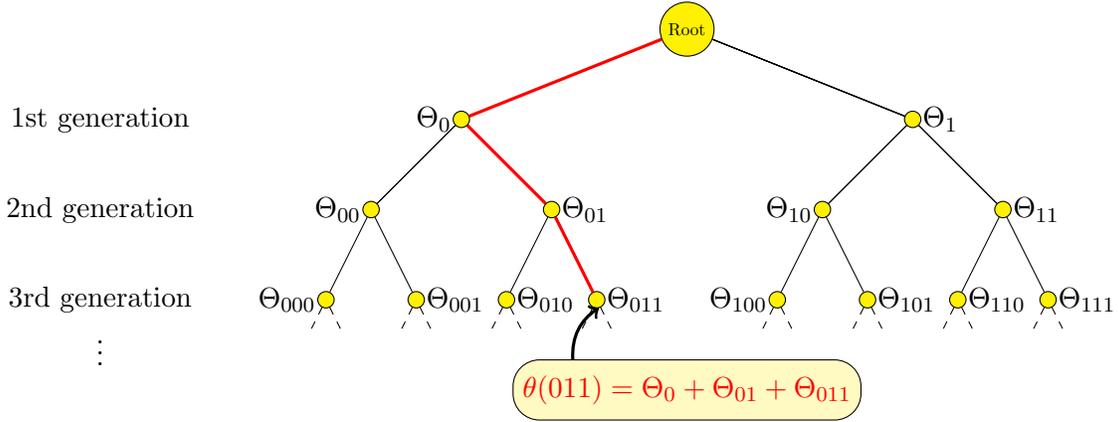
\begin{figure}[t] 
\centering
\begin{tikzpicture}[xscale=0.6,yscale=0.6] 
\tikzstyle{fleche}=[->,very thick,rounded corners=4pt];
\tikzstyle{sommet}=[circle,draw,fill=yellow,scale=0.6] 
\tikzstyle{noeud}=[minimum width=2cm, minimum height=0.8cm, rectangle, rounded corners=10pt,draw,fill=yellow!30,text=red,font=\bfseries] 
\node[sommet] (R) at (0,10){Root};
\node[sommet] (R0) at (-5,8){}; 
%\draw (R0)   node[left] {$\Theta_0$};
\node[sommet] (R1) at (5,8){};
\node[sommet] (R00) at (-7,6){};
%\draw (R00)   node[left] {$\Theta_0$};
\node[sommet] (R01) at (-3,6){};
\node[sommet] (R10) at (3,6){};
\node[sommet] (R11) at (7,6){};
\node[sommet] (R000) at (-8,4){};
\node[sommet] (R001) at (-6,4){};
\node[sommet] (R010) at (-4,4){};
\node[sommet] (R011) at (-2,4){};
\node[sommet] (R100) at (2,4){};
\node[sommet] (R101) at (4,4){};
\node[sommet] (R110) at (6,4){};
\node[sommet] (R111) at (8,4){};
\node (R0000) at (-8.5,3){};
\node  (R0001) at (-7.5,3){};
\node  (R0010) at (-6.5,3){};
\node  (R0011) at (-5.5,3){};
\node  (R0100) at (-4.5,3){};
\node  (R0101) at (-3.5,3){};
\node  (R0110) at (-2.5,3){};
\node  (R0111) at (-1.5,3){};
\node (R1000) at (1.5,3){};
\node (R1001) at (2.5,3){};
\node (R1010) at (3.5,3){};
\node (R1011) at (4.5,3){};
\node  (R1100) at (5.5,3){};
\node  (R1101) at (6.5,3){};
\node  (R1110) at (7.5,3){};
\node  (R1111) at (8.5,3){};
\node  (time1) at (-13,8){1st generation};
\node  (time2) at (-13,6){2nd generation};
\node  (time3) at (-13,4){3rd generation};
\node  (time4) at (-13,3){$\vdots$};
\draw (R) -- (R0) node[ left] {$\Theta_{0}$} -- (R00) node[ left] {$\Theta_{00}$} -- (R000) node[left] {$\Theta_{000}$};
\draw (R) -- (R0) -- (R00) -- (R001) node[right] {$\Theta_{001}$};
\draw (R) -- (R0) -- (R01) node[right] {$\Theta_{01}$} -- (R010) node[right] {$\Theta_{010}$};
\draw (R) -- (R0) -- (R01) -- (R011) node[right] {$\Theta_{011}$} ;
\draw (R) -- (R1) -- (R10) -- (R100) node[left] {$\Theta_{100}$};
\draw (R) -- (R1) node[right] {$\Theta_{1}$}-- (R10) node[left] {$\Theta_{10}$}-- (R101)node[right] {$\Theta_{101}$};
\draw (R) -- (R1) -- (R11) node[right] {$\Theta_{11}$}-- (R110)node[right] {$\Theta_{110}$};
\draw (R) -- (R1) -- (R11) -- (R111)node[right] {$\Theta_{111}$};
\draw[color=red,very thick]  (R) -- (R0) -- (R01) -- (R011);
\node[noeud] (theta) at (0,2) {$\theta(011)=\Theta_0+\Theta_{01}+\Theta_{011}$};
\draw[fleche] (theta) to[bend left] (R011);
\draw[dashed] (R000)--(R0000);
\draw[dashed] (R000)--(R0001);
\draw[dashed] (R001)--(R0010);
\draw[dashed] (R001)--(R0011);
\draw[dashed] (R010)--(R0100);
\draw[dashed] (R010)--(R0101);
\draw[dashed] (R011)--(R0110);
\draw[dashed] (R011)--(R0111);
\draw[dashed] (R100)--(R1000);
\draw[dashed] (R100)--(R1001);
\draw[dashed] (R101)--(R1010);
\draw[dashed] (R101)--(R1011);
\draw[dashed] (R110)--(R1100);
\draw[dashed] (R110)--(R1101);
\draw[dashed] (R111)--(R1110);
\draw[dashed] (R111)--(R1111);
\end{tikzpicture}
\caption{Representation of the branching random walk on the dyadic tree $\T$. In red, the genealogy of the particle $(011)$.}
\label{tree}
\end{figure}

We consider a family of independent complex Gaussian random variables $(\Theta_{u})_{u\in \mathbb{T}}   $    indexed by the nodes of this tree  and identically distributed with common law $\Theta$, whose real part is independent of the imaginary part  (see Figure \ref{tree}).

  If $u=(u_1, \cdots, u_n)\in  \mathbb{T}_{\leq n } $ and $l \leq n$  or $u\in\T$, we set $u_{|l}= (u_1, \cdots, u_l)\in  \mathbb{T}_{\leq l }$.  

 Each particle $u\in \mathbb{T}_{\leq n }$ is given a complex weight $\theta(u)$ corresponding to the sum of the weights $\Theta_e$ encountered along the shortest path joining $u$ to the root (see Figure \ref{tree}), i.e.
 $$\theta(u)=\sum_{l\leq n}\Theta_{u_{|l}}.$$ 
We define the real and imaginary part of $\theta$
$$V(u)={\rm Re}(\theta(u))\quad \text{ and }X(u)={\rm Im}(\theta(u)),$$
which are independent. We assume that the real part of the branching random walk is critical
\begin{equation}
\E\Big(\sum_{|z|=1}\ee^{-V(z)}\Big)=1,\qquad \E\Big(\sum_{|z|=1}V(z)\ee^{-V(z)}\Big)=0.
\end{equation}
Put in other words, the real part of $\Theta$ has law $\sqrt{2\ln 2}\mathcal{N}(0,1)+2\ln 2 $ where $\mathcal{N}(0,1)$ is a standard Gaussian random variable. The imaginary part of $\Theta$ has law $\mathcal{N}(0,1)$ and is independent of its real part.

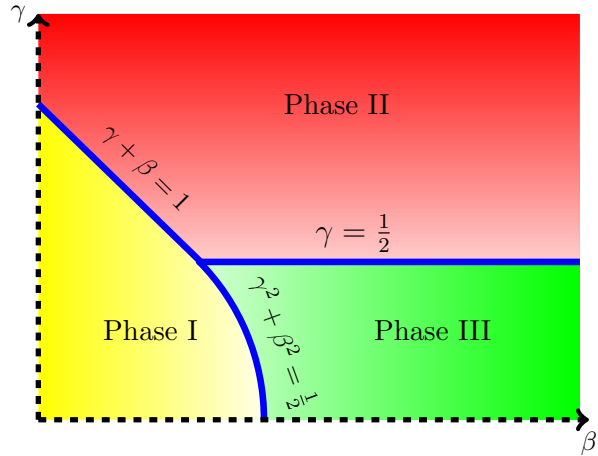
\begin{figure}[!h]
\centering
\begin{tikzpicture}[xscale=0.6,yscale=0.6] 
\shade[left color=green!10,right color=green] (3,0) -- (12,0) -- (12,6) -- (3,6) ;
\shade[top color=red, bottom color=red!20] (3.5,3.5) -- (12,3.5) -- (12,9) -- (0,9) -- (0,7);
\shade[left color=yellow,right color=yellow!10] (5,0) arc (0:45:5) -- (0,7) -- (0,0) -- (5,0);
%\draw[color=green!40,fill=green!40] (3,0) -- (12,0) -- (12,6) -- (3,6) ;
%\draw[color=red!40,fill=red!40] (3.5,3.5) -- (12,3.5) -- (12,8) -- (0,8) -- (0,7);
%\draw[color=yellow!40,fill=yellow!40] (5,0) arc (0:45:5) -- (0,7) -- (0,0) -- (5,0);
\draw [line width=2.5pt,color=blue](3.53,3.5) -- (12,3.5);
\draw [line width=2.5pt,color=blue](3.6,3.5) -- (0,7);
\draw [line width=2.5pt,color=blue](5,0)  arc (0:45:5);
\draw[style=dashed,line width=2pt,->] (0,0) -- (0,9) node[left]{{\small $\gamma$}};
\draw[style=dashed,line width=2pt,->] (0,0) -- (12.2,0)node[below]{{\small $\beta$}};
\draw (7,3.5)  node[above]{{\large $\gamma=\frac{1}{2}$}};
\draw (2,5.2)  node[above,rotate=-45,line width=2pt]{{\small $\gamma+\beta=1$}};
\draw (4.8,1.5)  node[above,rotate=-70,line width=2pt]{{\small $\gamma^2+\beta^2=\frac{1}{2}$}};
\draw (1,2)  node[right]{{  Phase I}};
\draw (5,7)  node[right]{{  Phase II}};
\draw (7,2)  node[right]{{  Phase III}};
\end{tikzpicture}
\caption{Phase diagram}
\label{diagram}
\end{figure}

\medskip
For any $(\gamma,\beta) \in \r_+^2$,  we define:
\begin{equation}\label{def:martingal}
M_n^{\gamma,\beta}:= \sum_{|z|=n}\ee^{-\gamma V(z)+ i\beta\sqrt{2\ln 2}X(z)}.
\end{equation}
One gets a martingale out of $(M_n^{\gamma,\beta})_n$ by renormalizing by the mean, namely by considering $(M_n^{\gamma,\beta}/\E[M_n^{\gamma,\beta}])_n$. It is natural to wonder for which values  of the parameters $(\gamma,\beta)$ this renormalization by the mean gives a martingale converging almost surely towards a non trivial limit $M^{\gamma,\beta}$.  The real case $\beta=0$ has given rise to an extensive literature ranging  from the study of Mandelbrot's multiplicative cascades (see \cite{mandelbrot,KP,Barral1,BKNSW} among many others) to some extensions like the study of Branching random walks (see \cite{biggins,Biki,AS,HuShi} and references therein). The complex case has been answered in \cite{bar:comp1,bar:comp2,lrv} (actually, \cite{lrv} is concerned with Gaussian multiplicative chaos, which is a different context, but the methods apply to our context). The reader may have in mind the resulting phase diagram in Figure \ref{diagram}. 
It is proved in  \cite{bar:comp1,bar:comp2} that  we have almost sure convergence towards a non trivial limit in the so-called phase I, i.e. 
 \begin{equation}\label{phase1}
 \Big(|\gamma|\leq 1/2\text{ and }\gamma^2+\beta^2<1/2\Big)\text{ or }\Big(\gamma\in ]\frac{1}{2},1[\text{ and }\gamma+\beta<1\Big)
\end{equation}
 and the boundary case 
 \begin{equation}\label{boundary}
\gamma\in ]\frac{1}{2},1[\text{ and }\gamma+\beta=1
\end{equation}
is treated in \cite{lrv} on a related model, Gaussian multiplicative chaos. In all other cases, a renormalization by the mean is a martingale that does not converge to something non trivial.

A mathematical understanding of the limiting objects we get is therefore a natural question. In particular, one may look at \eqref{def:martingal} more generally as a complex measure on $[0,1]$. The canonical way to do so is to use the dyadic decomposition of the reals that belong to $[0,1]$. More specifically, if $u \in \mathbb{T}_{\leq n}$ then we set $t_u=\sum_{i=1}^n \frac{u_i}{2^i}\in [0,1]$; the complex measure $M_n^{\gamma,\beta}(dx)$ is then defined on $[0,1]$ by       
\begin{equation*}
M_n^{\gamma,\beta}(dx)= \sum_{|z|=n}\ee^{-\gamma V(z)+ i\beta\sqrt{2\ln 2}X(z)}\mathbf{1}_{]t_z,t_z+2^{-n}]}(x)\,2^n\,dx.
\end{equation*} 
In the case when $\beta=0$, these complex measures turn out to be   positive random measures and an important question is to study the  diffusivity properties of the scaling limit. It is  known nowadays that the limiting measures are diffuse in phase I \cite{KP}  or its boundary (i.e. $\gamma=1,\beta=0$, see \cite{BKNSW}): they are atom free. In the complex case (i.e. $\beta \not = 0$), this question makes sense when asking whether the family of random functions
$$t\in[0,1]\mapsto M_n^{\gamma,\beta}[0,t]$$
converge uniformly towards a continuous limit after renormalization by the mean of the total mass $\E(M_n^{\gamma,\beta}[0,1])$.

In the inner phase I,  it is proved in  \cite{bar:comp1,bar:comp2} that we have uniform convergence of this family of functions (even in the H\"older sense) as an application of the Kolmogorov criterion. In this paper, we aim at answering this question in the boundary case \eqref{boundary}, i.e. the frontier of phases I/II excluding the extremal points. In this case, the Kolmogorov criterion (or refined versions) breaks down.  Hence establishing uniform convergence is more difficult as one can no longer rely on general machinery on convergence or tightness of stochastic processes. Let us also mention that we have chosen a normalization of our parameters $(\gamma,\beta)$ so that the mean in the boundary case is exactly one, that is
\begin{Proposition}
For any $\gamma,\, \beta \in (\frac{1}{2},1)$ with $\beta +\gamma=1$ and any Borel set $A$ of $[0,1]$, the family $(M_n^{\gamma,\beta}(A))_n$  is a complex valued martingale with mean $1$.
\end{Proposition}

Now, we can state the main results of this paper:

\begin{theorem}\label{th:main}
Let $(\gamma,\beta)$ belong to the boundary of phases I/II, i.e. \eqref{boundary}. The sequence of functions $t \mapsto M_n^{\gamma,\beta}[0,t]$ converges almost surely in the space of continuous functions towards some random continuous function $M^{\gamma,\beta}$. There exists some (non explicit) constant $\xi>0$ (given by $\gamma \frac{1-\eta}{2}$ where $\eta$ appears in proposition \ref{boundness}) such that for all $\epsilon>0$ the function $M^{\gamma,\beta}$ satisfies almost surely the following modulus estimate   
\begin{equation}\label{eq:mainmodulus} 
|M^{\gamma,\beta}[s,t]| \leq C \frac{1}{ ( \ln (1+ \frac{1}{|t-s|}) )^{\xi-\epsilon} }, \quad s,t \in [0,1]
\end{equation}
where $C>0$ is some random constant.
\end{theorem}
\begin{Remark}
We could in principle give an explicit formula for $\xi$ (as a function of $\gamma$) but since our method is not optimal, we did not try to keep track of this information.   
\end{Remark}
 
 As mentioned in the introduction, it turns out that the space of continuous functions is the right space to study convergence of $ M_n^{\gamma,\beta}$. Indeed, though for each $n \geq 1$  $ M_n^{\gamma,\beta}$ is a complex measure, one can show that the limiting function $M^{\gamma,\beta}$ is not of finite variation (see subsection \ref{infinitevar})
\begin{corollary}\label{coro:main}
Let $(\gamma,\beta)$ belong to the boundary of phases I/II, i.e. \eqref{boundary}. The mapping $t \mapsto M^{\gamma,\beta}[0,t]$ is not of finite variation almost surely.
\end{corollary}

\subsection{Related models and open problems}
%%%%%%%%%%%%%%%%%%%%%%%%
 The authors of \cite{lrv} studied complex Gaussian multiplicative chaos (GMC) in dimension $d$, i.e. the renormalization theory of $e^{ \sqrt{2d}\gamma X(x) + i \sqrt{2d} \beta Y(x) } dx$ where $dx$ is the Lebesgue measure, $\gamma, \beta \in \r_+^2$, $X$ and $Y$ are two independent log-correlated fields on some domain $D \subset \r^d$ (of course, we could absorb the $\sqrt{2d}$ in the parameters $\gamma$ and $\beta$ but we will not do so in order to compare complex GMC with the complex branching random walk studied in this paper). In particular, the work \cite{lrv}  enables to define $e^{ \sqrt{2}\gamma X(x) + i \sqrt{2} \beta Y(x) } dx$ where $X$ and $Y$ are two independent log-correlated fields on $\r$ with covariance given by $\E[X(x)X(y)]= \ln \frac{1}{|y-x|}$ and when $(\gamma, \beta)$ satisfies the condition 
$$\gamma+\beta=1, \; \gamma \in (\frac{1}{2},1) \quad \text{I/II}.$$ 
Condition I/II refers to the so-called frontier between phase I and phase II: see \cite{lrv}. More precisely, if $X_\varepsilon, Y_\varepsilon$ are appropriate cut-off approximations of $X,Y$ (as $\varepsilon$ goes to $0$) then, if $(\gamma,\beta)$ are in phase I/II, the random distribution
\begin{equation*}
\bar{M}_\varepsilon^{\gamma,\beta}:=e^{ \sqrt{2}\gamma X_\varepsilon(x) + i \sqrt{2} \beta Y_\varepsilon (x) -(\gamma^2-\beta^2) \E[X_\epsilon(x)^2]} dx
\end{equation*}
converges almost surely (in the space of distributions) towards a distribution $\bar{M}^{\gamma,\beta}$. In this case, the operator $\bar{M}^{\gamma,\beta}$ is a conformally invariant boundary operator in the framework of $2D$ string theory on the upper half plane (see \cite{lrv}). The work \cite{lrv} did not establish but conjectured that convergence holds in the space of continuous functions. %We will reinforce this conjecture below by adding an estimate on the modulus of the limiting object $\bar{M}^{\gamma,\beta}$.  

Now, one expects that the complex branching random walk of this paper and complex GMC have a similar behaviour. One justification for this is that the branching random walk also has logarithmic correlations but with respect to the underlying ultrametric distance on the tree $\mathbb{T}$. Hence, theorem \ref{th:main} gives additional support to the conjecture of \cite{lrv}. 

\medskip 

Another important question on this topic is the following. When the renormalized martingale defined by \eqref{def:martingal} converges towards a non trivial limit, it is readily seen that the limit satisfies a distributional equation of the type
$$M^{\gamma,\beta}=\sum_{|z|=1}T(z)M^{\gamma,\beta}(z)$$
where $(T(z))_{|z|=1}$ are complex random variables independent of  $(M^{\gamma,\beta}(z))_{|z|=1}$, which are i.i.d. random variables with law $M^{\gamma,\beta}$. Such an equation is known under the name of smoothing transform and has been extensively studied in the case when $(T(z))_{|z|=1}$ are positive  (see \cite{Biggf,durrett,Kyp}), real valued \cite{meiners}. We would like to emphasize the fact that understanding this equation in the complex case is an important point. Let us mention some work to appear \cite{meiners2} in this perspective.
\subsection*{Acknowledgements}
The authors would like to thank Matthias Meiners for interesting discussions.   

%%%%%%%%%%%%%%%%%%%%%%%%%%%%%%%%%%%%%%%%%%%%%%%%%%%%%%%%%%%%%%%%%%%%%%%%%%%%%%%%
\section{Proof of the main result}
%%%%%%%%%%%%%%%%%%%%%%%%%%%%%%%%%%%%%%%%%%%%%%%%%%%%%%%%%%%%%%%%%%%%%%%%%%%%%%%%

 In this section, we first introduce some further notations and then give the proof of Theorem \ref{th:main} based on some auxiliary results, the proofs of which are postponed to the following sections in the paper.
 
 \subsection{Convention on notations}
 In the sequel, we will denote by $c$ a generic positive constant which can change from line to line; we will also denote $\eta$ a generic constant belonging to $(0,1)$ whose value can also change from line to line. Usually, it will be clear from the context that $c$ or $\eta$ can depend on other constants but to lighten notations, we will make this dependence implicit.

 \subsection{Notations}
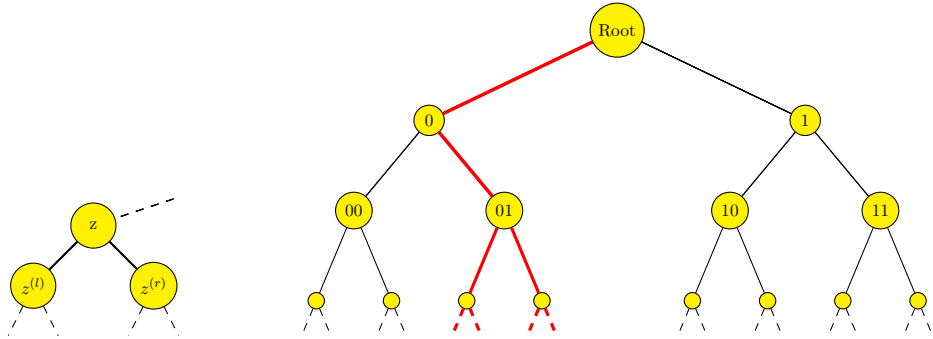
\begin{figure}[t] 
\centering 
\begin{tikzpicture}[xscale=0.4,yscale=0.4] 
\tikzstyle{fleche}=[->,very thick,rounded corners=3pt];
\tikzstyle{sommet}=[circle,draw,fill=yellow,scale=0.6,minimum height=1cm] 
\tikzstyle{noeud}=[minimum width=2cm, minimum height=0.8cm, rectangle, rounded corners=10pt,draw,fill=yellow!30,text=red,font=\bfseries] 
\node (R) at (3,9){};
\node[sommet] (R1) at (0,8){z};
\node[sommet] (R10) at (-2,6){$z^{(l)}$};
\node[sommet] (R11) at (2,6){$z^{(r)}$};
\node  (R100) at (-3,4){};
\node  (R101) at (-1,4){};
\node  (R110) at (1,4){};
\node  (R111) at (3,4){};
\draw [dashed] (R) -- (R1) -- (R10) -- (R100);
\draw  [dashed] (R) -- (R1) -- (R10) -- (R101);
\draw[thick] (R1) -- (R10);
\draw[thick] (R1) -- (R11);
\draw [dashed] (R) -- (R1) -- (R11) -- (R110);
\draw [dashed] (R) -- (R1) -- (R11)-- (R111);
\end{tikzpicture}
%\caption{Node and children.}
%\label{children}
%
\hspace{1cm}
\begin{tikzpicture}[xscale=0.5,yscale=0.6] 
\tikzstyle{fleche}=[->,very thick,rounded corners=4pt];
\tikzstyle{sommet}=[circle,draw,fill=yellow,scale=0.6] 
\tikzstyle{noeud}=[minimum width=2cm, minimum height=0.8cm, rectangle, rounded corners=10pt,draw,fill=yellow!30,text=red,font=\bfseries] 
\node[sommet] (R) at (0,10){Root};
\node[sommet] (R0) at (-5,8){0};
\node[sommet] (R1) at (5,8){1};
\node[sommet] (R00) at (-7,6){00};
\node[sommet] (R01) at (-3,6){01};
\node[sommet] (R10) at (3,6){10};
\node[sommet] (R11) at (7,6){11};
\node[sommet] (R000) at (-8,4){};
\node[sommet] (R001) at (-6,4){};
\node[sommet] (R010) at (-4,4){};
\node[sommet] (R011) at (-2,4){};
\node[sommet] (R100) at (2,4){};
\node[sommet] (R101) at (4,4){};
\node[sommet] (R110) at (6,4){};
\node[sommet] (R111) at (8,4){};
\node (R0000) at (-8.5,3){};
\node  (R0001) at (-7.5,3){};
\node  (R0010) at (-6.5,3){};
\node  (R0011) at (-5.5,3){};
\node  (R0100) at (-4.5,3){};
\node  (R0101) at (-3.5,3){};
\node  (R0110) at (-2.5,3){};
\node  (R0111) at (-1.5,3){};
\node (R1000) at (1.5,3){};
\node (R1001) at (2.5,3){};
\node (R1010) at (3.5,3){};
\node (R1011) at (4.5,3){};
\node  (R1100) at (5.5,3){};
\node  (R1101) at (6.5,3){};
\node  (R1110) at (7.5,3){};
\node  (R1111) at (8.5,3){};
\draw (R) -- (R0) -- (R00)  -- (R000);
\draw (R) -- (R0) -- (R00) -- (R001);
\draw (R) -- (R0) -- (R01) -- (R010);
\draw (R) -- (R0) -- (R01) -- (R011);
\draw (R) -- (R1) -- (R10) -- (R100);
\draw (R) -- (R1) -- (R10)-- (R101);
\draw (R) -- (R1) -- (R11) -- (R110);
\draw (R) -- (R1) -- (R11) -- (R111);
\draw[color=red,very thick]  (R) -- (R0) -- (R01) -- (R011);
\draw[color=red,very thick]  (R) -- (R0) -- (R01) -- (R010);
\draw[dashed] (R000)--(R0000);
\draw[dashed] (R000)--(R0001);
\draw[dashed] (R001)--(R0010);
\draw[dashed] (R001)--(R0011);
\draw[dashed,color=red,very thick] (R010)--(R0100);
\draw[dashed,color=red,very thick] (R010)--(R0101);
\draw[dashed,color=red,very thick] (R011)--(R0110);
\draw[dashed,color=red,very thick] (R011)--(R0111);
\draw[dashed] (R100)--(R1000);
\draw[dashed] (R100)--(R1001);
\draw[dashed] (R101)--(R1010);
\draw[dashed] (R101)--(R1011);
\draw[dashed] (R110)--(R1100);
\draw[dashed] (R110)--(R1101);
\draw[dashed] (R111)--(R1110);
\draw[dashed] (R111)--(R1111);
\end{tikzpicture}
\caption{Left: node and children. Right: subtree rooted  at $(01)$ in red.}
\label{rootened}
\end{figure}

%%%%%%%%%%%%%%%%%%%%
For any $z\in \mathbb{T}= \bigcup_n\mathbb{T}_{\leq n}$, we let $z^{(l)}$ be the left child of $z$ and $z^{(r)}$ be the right child of $z$ (see Figure \ref{rootened}). Notice that if $|z|=k$ then $|z^{(l)}|=|z^{(r)}|=k+1$. Furthermore, to simplify the notations, we mention that we will implicitly assume in the following that the notation $u_{|k}^{(l)}$ means $(u_{|k})^{(l)}$ (note that it is not clear otherwise whether one must take the child  after or before restricting to the k-th node of $u$). The same convention holds for the right child.

For any $u\in \mathbb{T}_{\leq n}$, we further consider the sub-tree $\mathbb{T}_u$ of $\mathbb{T}$ rootened in $u$ (see Figure \ref{rootened}), namely
\begin{equation}
\mathbb{T}_u:=\{ v\in \mathbb{T},\, v_{|n}=u\},
\end{equation}
 Finally, we define an order on $\T$. We write $u \geq z$  if $u$ is a descendant of $z$ in the tree $\mathbb{T}$, i.e. if $|u| \geq |z|$ and $u_{| |z|}=z$ (we adopt similar conventions for $\leq$). 
  
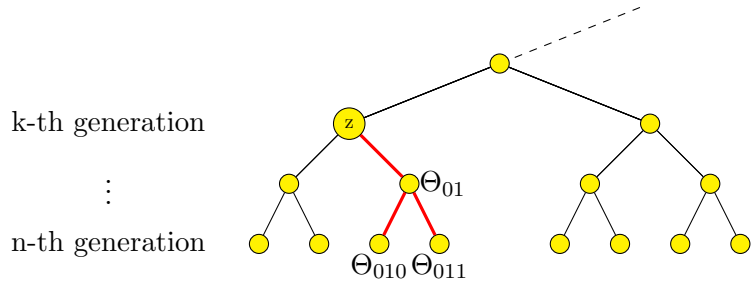
\begin{figure}[!h] 
\centering
\begin{tikzpicture}[xscale=0.4,yscale=0.4] 
\tikzstyle{fleche}=[->,very thick,rounded corners=4pt];
\tikzstyle{sommet}=[circle,draw,fill=yellow,scale=0.7] 
\tikzstyle{noeud}=[minimum width=2cm, minimum height=0.8cm, rectangle, rounded corners=10pt,draw,fill=yellow!30,text=red,font=\bfseries] 
\node (1R) at (5,12){};
\node[sommet] (R) at (0,10){};
\node[sommet] (R0) at (-5,8){z};
\node[sommet] (R1) at (5,8){};
\node[sommet] (R00) at (-7,6){};
\node[sommet] (R01) at (-3,6){};
\node[sommet] (R10) at (3,6){};
\node[sommet] (R11) at (7,6){};
\node[sommet] (R000) at (-8,4){};
\node[sommet] (R001) at (-6,4){};
\node[sommet] (R010) at (-4,4){};
\node[sommet] (R011) at (-2,4){};
\node[sommet] (R100) at (2,4){};
\node[sommet] (R101) at (4,4){};
\node[sommet] (R110) at (6,4){};
\node[sommet] (R111) at (8,4){};
\node  (time1) at (-13,8){k-th generation};
\node  (time2) at (-13,6){$\vdots$};
\node  (time3) at (-13,4){n-th generation};
\draw (R) -- (R0) -- (R00) -- (R000);
\draw (R) -- (R0) -- (R00) -- (R001) ;
\draw (R) -- (R0) -- (R01) node[right] {$\Theta_{01}$} -- (R010) node[below] {$\Theta_{010}$};
\draw (R) -- (R0) -- (R01) -- (R011) node[below] {$\Theta_{011}$} ;
\draw (R) -- (R1) -- (R10) -- (R100) ;
\draw (R) -- (R1) -- (R10) -- (R101);
\draw (R) -- (R1) -- (R11) -- (R110);
\draw (R) -- (R1) -- (R11) -- (R111);
\draw[color=red,very thick]    (R0) -- (R01) -- (R011);
\draw[color=red,very thick]     (R01) -- (R010);
\draw[dashed] (1R)--(R);
\end{tikzpicture}
\caption{In red, representation of $M_n^{\gamma,\beta}(z)$: we sum over the red paths rooted at $z$ and going down up to the $n$-th generation. The mass of each path is the sum over the variables $\Theta$ encountered along the way down.}
\label{Mnz}
\end{figure}

We set $\mathbb{M}_n:= \inf_{|u|=n}V(u)$. For any $n \in \N$ and $u\in \mathbb{T}_{\leq l}$ with $l \leq n$, we consider the mass of the subtree rootened at $z$ up to generation $n$ (see Figure \ref{Mnz})
\begin{align}\label{defMnz}
M_n^{\gamma,\beta}(u):=&\ee^{\gamma V(u)- i\beta\sqrt{2\ln 2}X(u)}M_n^{\gamma,\beta}[t_u,t_u+\frac{1}{2^{l}}]\\
= &\sum_{|z|=n,\, z_{|l}=u}  \ee^{-\gamma [V(z)-V(u)]+ i\beta\sqrt{2\ln 2}[X(z)-X(u)]}.\nonumber
\end{align}

\begin{figure}[!h] 
\centering
\begin{tikzpicture}[xscale=0.4,yscale=0.4] 
\tikzstyle{fleche}=[->,very thick,rounded corners=4pt];
\tikzstyle{sommet}=[circle,draw,fill=yellow,scale=0.7] 
\tikzstyle{noeud}=[minimum width=2cm, minimum height=0.8cm, rectangle, rounded corners=10pt,draw,fill=yellow!30,text=red,font=\bfseries] 
\node (1R) at (5,12){};
\node[sommet] (R) at (0,10){};
\node[sommet] (R0) at (-5,8){};
\node[sommet] (R1) at (5,8){};
\node[sommet] (R00) at (-7,6){};
\node[sommet] (R01) at (-3,6){};
\node[sommet] (R10) at (3,6){};
\node[sommet] (R11) at (7,6){};
\node[sommet] (R000) at (-8,4){};
\node[sommet] (R001) at (-6,4){};
\node[sommet] (R010) at (-4,4){};
\node[sommet] (R011) at (-2,4){};
\node[sommet] (R100) at (2,4){};
\node[sommet] (R101) at (4,4){u};
\node[sommet] (R110) at (6,4){};
\node[sommet] (R111) at (8,4){};
\node  (time1) at (-13,8){k-th generation};
\node  (time2) at (-13,6){$\vdots$};
\node  (time3) at (-13,4){n-th generation};
\draw (R) -- (R0) -- (R00) -- (R000);
\draw (R) -- (R0) -- (R00) -- (R001) ;
\draw (R) -- (R0) -- (R01)  -- (R010)  ;
\draw (R) -- (R0) -- (R01) -- (R011)   ;
\draw (R) -- (R1) -- (R10) -- (R100) ;
\draw (R) -- (R1) -- (R10) -- (R101);
\draw (R) -- (R1) -- (R11) -- (R110);
\draw (R) -- (R1) -- (R11) -- (R111);
\draw[color=red,very thick]    (R) -- (R1) -- (R10) -- (R101);
\draw[dashed,color=red,very thick] (1R)--(R);
%bulles
\draw[color=blue,thick, rounded corners=8pt] (R) -- (-6,8.5) -- (-9,3.5) -- (-1,3.5) -- (-3.5,8)--   (R);
\draw[color=blue,thick, rounded corners=8pt] (R10) -- (1,3.5) -- (2.5,3.5) --  (R10);
\end{tikzpicture}
\caption{The decomposition \eqref{decompMn} amounts to summing over all the subtrees (wrapped in blue) located on the left of the red path.}
\label{lefttree}
\end{figure}
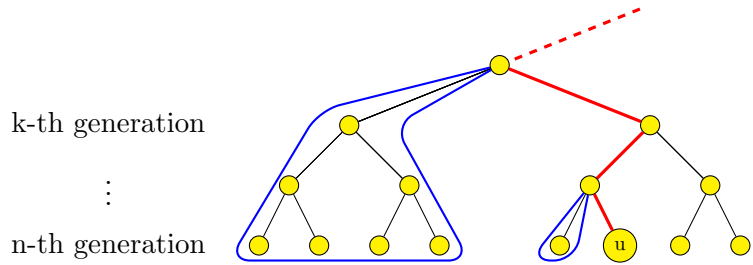

With  these notations, we have for $|u|=n$
\begin{equation}\label{decompMn}
M_n^{\gamma,\beta}[0,t_u]= \sum_{k=0}^{n-1} \ee^{-\gamma V(u_{|k}^{(l)}) +i\beta\sqrt{2\ln 2}X(u_{|k}^{(l)})} M_{n}^{\gamma,\beta}(u_{|k}^{(l)}) 1_{\{ u_{| k+1}\neq u_{|k}^{(l)}  \}}.
\end{equation} 
This can be seen by summing over all the subtrees located on the left-hand side of the path joining the root to the particle $u$ (see Figure \ref{lefttree}).

In fact, we can extend this decomposition: if $l \leq n$ and $|u|=n$  
\begin{align*}
M_n^{\gamma,\beta}[t_{u_{|l}},t_u]=& \ee^{-\gamma V(u_{|l}) +i\beta\sqrt{2\ln 2}X(u_{|l})} \\
&\times\sum_{k=0}^{n-l-1} \ee^{-\gamma (V(u_{|k+l}^{(l)})-V(u_{|l})   )  +i\beta\sqrt{2\ln 2}(  X(  u_{|k+l}^{(l)} )-X(u_{|l}))} M_{n}^{\gamma,\beta}(u_{|k+l}^{(l)}) 1_{\{ u_{|k+l+1}\neq u_{|k+l}^{(l)}  \}}  .
\end{align*}
 
Finally, we introduce the following quantity for $|u|=l$, $l \leq n$ and $p \geq 0$
\begin{equation}\label{defquantity}
||  M^{\gamma,\beta}_{n,p}(u) ||_{\infty} = \max_{|z|=n, z_{|l}=u}  \sum_{k=0}^{n-l} \ee^{-\gamma (V(z_{|k+l}^{(l)})-V(u))} |M_{n+p}^{\gamma,\beta}(z_{|k}^{(l)})| 1_{\{ z_{|k+l+1}\neq z_{|k+l}^{(l)}  \}} + \ee^{-\gamma (V(z)-V(u) )}  | M^{\gamma,\beta}_{n+p}(z)  |
\end{equation}
We also set $||  M^{\gamma,\beta}_{n}(u)||_{\infty}  :=  ||M^{\gamma,\beta}_{n,p=0}(u)||_{\infty}  $ and $||M^{\gamma,\beta}_{n,p}||_{\infty}:= ||M^{\gamma,\beta}_{n,p}(\emptyset)||_{\infty} $ .

 Note that by the triangle inequality we have
 \begin{equation}\label{comparaisondebut}
\sup_{|z|=n}|M_{n+p}^{\gamma,\beta}[0,t_z+\frac{1}{2^n}]| \leq ||  M^{\gamma,\beta}_{n,p} ||_{\infty} 
 \end{equation}
hence upper bounds on $||  M^{\gamma,\beta}_{n,p} ||_{\infty}$ lead to upper bounds on $ \sup_{|z|=n}|M_{n+p}^{\gamma,\beta}[0,t_z+\frac{1}{2^n}]|$ (which of course is equal to $\sup_{k \leq 2^n}|M_{n+p}^{\gamma,\beta}[0,\frac{k}{2^n}]|$). Finally, by the recursive structure on the tree, i.e. the subtree starting from any vertex $u$ has same distribution as the original tree, the variable $||  M^{\gamma,\beta}_{n,p}(u) ||_{\infty} $ has the same distribution as $||M^{\gamma,\beta}_{r,p}||_{\infty}$ for $r=n-l$.

\subsection{Many to one lemma and useful estimates}

We introduce the centered standard Gaussian walk $(S_n)_{n \geq 1}$. When the walk starts from a point $x$, we denote the associated probability measure $\P_x$ and the corresponding measure $\E_x$. When $x=0$, we will omit the subscript. We will denote $\underline{S}_n$ the infimum of the walk on the set $\lbrace 1, \ldots, n \rbrace$.

Finally, we recall the many to one lemma which is very useful in the context of branching random walks; for all function $F$
\begin{equation}\label{manytoone}
\E\Big(\sum_{|z|=n}  F( x+ V(z))\ee^{-V(z)}\Big)= \E_x   [F(S_n)]    
\end{equation}

Recall also that it is proven in \cite{AS} that:

-(see Lemma 2.4) there exists $c>0$ such that for any $n\in\N$, $x\geq 0$ and $a\leq b\in \r$, 
\begin{equation}\label{estimAid}
\P_x\left( \min_{j\leq n}{S}_n\geq 0,\, S_n\in [a,b]\right)\leq \frac{c(1+x)(1+b)(1+(b-a))}{n^{\frac{3}{2}}}.
\end{equation}

-(see Lemma B.2) for any $\kappa>0$ there exists $c(\kappa)>0$ such that for any $x\geq 0$ 
\begin{equation}
\label{estimAid2}
\E_x\left[\sum_{l\geq 0} \ee^{-\kappa S_l} 1_{\{ \min_{j\leq l}S_j\geq 0  \}}\right] \leq c(\kappa)
\end{equation}
In fact, lemma 2.4 and Lemma B.2 of \cite{AS} are much more general and concerns more general walks than the Gaussian one.

\subsection{Proof of Theorem \ref{th:main}}
%%%%%%%%%%%%%%%%%%%%%%%%%

\subsubsection{Idea of the proof}

We first explain the main idea behind the proof of theorem \ref{th:main}. Essentially, the proof relies on the fact that one can extract from a uniform (in $n$) bound on the supremum of $t \mapsto M^{\gamma,\beta}_n[0,t]$ a uniform (in $n$) bound on the increments of $t \mapsto M^{\gamma,\beta}_n[0,t]$. More precisely, the main estimate of this paper is given by proposition \ref{boundness} below which gives a uniform bound on $||  M^{\gamma,\beta}_{n,p} ||_{\infty}$ for all $n,p$. Recall that, by \eqref{comparaisondebut}, one should see inequality \eqref{eqboundness} below as an estimate on  the supremum of $t \mapsto M^{\gamma,\beta}_{n+p}[0,t]$. In fact, we will only use the estimate for $p=0$. Now, by the recursive structure on the tree, i.e. the subtree starting from any vertex $u$ has same distribution as the original tree, one can transfer such global estimates to estimates on the modulus of continuity of $M_{n+p}^{\gamma,\beta}$: this is essentially the content of inequality \eqref{thebound} (recall that $ ||M^{\gamma,\beta}_{r}(u)||_{\infty}$ is defined in \eqref{defquantity} and has same distribution as $||M^{\gamma,\beta}_{n}||_{\infty}$ with $n=r-l$ when $|u|=l$). This is sufficient to get tightness in the space of continuous functions and then almost sure convergence is a consequence of a non trivial theorem on Banach valued martingales in \cite{kruk} (the main theorem of  \cite{kruk} states that an $L^1$ bounded Banach valued martingale converges almost surely if and only if it is tight: here we will work in the Banach space of continuous functions).

\subsubsection{The proof}

Now, let us fix $l \geq 1$ and $r$ such that $l \leq r$.  % Finally, for $|u|=l$ we set for all $r \geq l$: 
%%\begin{equation*}
%%||  M_r(u) ||_{\infty}= \max_{|z|=r, z_{|l}=u} \left |\sum_{k=0}^{r-l} \ee^{-\gamma (V(z_{k+l}^{(l)})-V(u))+i\beta\sqrt{2\ln 2}(X(z_{k+l}^{(l)})-X(u))} M_{r}^{\gamma,\beta}(z_{k}^{(l)}) 1_{\{ z_{k+1}\neq z_k^{(l)}  \}} + \ee^{-\gamma V(z)}  M^{\gamma,\beta}_{r}(z)  \right |
%%\end{equation*}
%%Note that we have
%\begin{equation*}
%||  M_r(u) ||_{\infty} = \max_{|z|=r, z_{|l}=u}  \sum_{k=0}^{r-l} \ee^{-\gamma (V(z_{k+l}^{(l)})-V(u))} |M_{r}^{\gamma,\beta}(z_{k}^{(l)})| 1_{\{ z_{k+1}\neq z_k^{(l)}  \}} + \ee^{-\gamma V(z)}  | M^{\gamma,\beta}_{r}(z)  |
%\end{equation*}
 Since
\begin{equation*}
 \max_{|z|=r, z_{|l}=u} \left |\sum_{k=0}^{r-l} \ee^{-\gamma (V(z_{k+l}^{(l)})-V(u))+i\beta\sqrt{2\ln 2}(X(z_{k+l}^{(l)})-X(u))} M_{r}^{\gamma,\beta}(z_{k}^{(l)}) 1_{\{ z_{|k+1}\neq z_{|k}^{(l)}  \}} + \ee^{-\gamma V(z)}  M^{\gamma,\beta}_{r}(z)  \right | \leq ||  M^{\gamma,\beta}_r(u) ||_{\infty},
\end{equation*}
for $|u|=l$ we have:
\begin{equation*}
\sup_{s,t \in [t_u, t_u+\frac{1}{2^l}]} |M_r^{\gamma,\beta}[s,t]|  \leq 2  e^{-\gamma V(u)} ||  M^{\gamma,\beta}_r(u) ||_{\infty}.
\end{equation*}
Also, for $|u|=l$ and $|u'|=l$ the right neighboor on the tree (i.e. $t_{u'}=t_u+\frac{1}{2^l}$), we have:
\begin{equation*}
\sup_{s,t \in [t_u, t_u+\frac{2}{2^l}]} |M_r^{\gamma,\beta}[s,t]|  \leq 2  e^{-\gamma V(u)} ||  M^{\gamma,\beta}_r(u) ||_{\infty}+   e^{-\gamma V(u')} ||  M^{\gamma,\beta}_r(u') ||_{\infty}.
\end{equation*}
If $|t-s| \leq \frac{1}{2^l}$ then there are two possibilities: 

\vspace{0.2 cm}

\noindent \emph{First case}: $s$ and $t$ lie in the same dyadic interval of the form $[\frac{k}{2^l}, \frac{k+1}{2^l}]$ in which case we have $|M_r^{\gamma,\beta}[s,t]| \leq |M_r^{\gamma,\beta}[\frac{k}{2^l},s]|+ |M_r^{\gamma,\beta}[\frac{k}{2^l},t]| $ 

\vspace{0.2 cm}

\noindent \emph{Second case}: $s$ and $t$ lies in some dyadic interval of the form $[\frac{k}{2^l}, \frac{k+1}{2^l}]$ and $t$ in the dyadic interval  $[\frac{k+1}{2^l}, \frac{k+2}{2^l}]$ in which case we have $|M_r^{\gamma,\beta}[s,t]| \leq |M_r^{\gamma,\beta}[\frac{k}{2^l},s]|+ |M_r^{\gamma,\beta}[\frac{k}{2^l},\frac{k+1}{2^l}]|  + |M_r^{\gamma,\beta}[\frac{k+1}{2^l},t]| $ 

\vspace{0.2 cm}

Therefore, we get the bound:

\begin{equation}\label{thebound}
 \sup_{|t-s|\leq \frac{1}{2^l}} |M_r^{\gamma,\beta}[s,t]| \leq 3  \sup_{|u|=l} e^{-\gamma V(u)} ||  M^{\gamma,\beta}_r(u) ||_{\infty} 
\end{equation}

Recall that the recursive structure of the tree entails that $ ||  M^{\gamma,\beta}_r(u) ||_{\infty} $ above is distributed for each $u$ like $ ||  M^{\gamma,\beta}_n ||_{\infty}$ with $n=r-l$.

Let $\epsilon>0$. We set $\delta_{l,\epsilon}= \frac{1}{l^{(1-\eta)/2-\epsilon} }$ where $\eta \in (0,1)$ is such that the bound \eqref{eqboundness} page \pageref{eqboundness} holds.
Finally, we introduce
\begin{equation*}
\mathbb{A}_{l,\epsilon}: \{ \forall |u|=l, \, \min_{j\leq l}V(u_{|j}) \geq -\log l   \}  \cap  \{ \frac{3c}{\delta_{l,\epsilon}}\sum_{|u|=l} (1+V(u)^\eta 1_{\min_{j\leq l}V(u_{|j}) \geq -\log l }) e^{-V(u)} \leq \frac{1}{l^{\epsilon/2}}   \} 
\end{equation*}
Recall that for all $\alpha \in (0,1)$, there exists $c>0$ such that (see \cite{AS} Lemma 2.3 for example)
\begin{equation*}
\E[ (  l^{\frac{1-\eta}{2}} \sum_{|u|=l} (1+V(u)^\eta 1_{\min_{j\leq l}V(u_{|j}) \geq  -\log l }) e^{-V(u)}  )    ]= \E\left[ l^{\frac{1}{2}}  \frac{(1+S_l^\eta)}{l^{\eta \over 2}} 1_{ \min_{j\leq l} S_j\geq -\log l  }    \right] \leq c(\log l).
\end{equation*} 
Therefore, we have for all
\begin{align*}
& \P(  \frac{1}{\delta_{l,\epsilon}}\sum_{|u|=l} (1+V(u)^\eta 1_{\min_{j\leq l} V(u_{|j}) \geq -\log l}) e^{-V(u)} \geq \frac{1}{l^{\epsilon/2}}   ) \\
& \leq    l^{\frac{1-\eta}{2}  -\frac{\epsilon}{2}  } \E[   \sum_{|u|=l} (1+V(u)^\eta 1_{\min_{j\leq u}V(u_{|j}) \geq   -\log l}) e^{-V(u)}      ]  \\
& \leq \frac{c\log l }{l^{\frac{\epsilon}{2}} }  .
\end{align*}
Then it leads to the estimate   
 \begin{equation*}
 \P( \mathbb{A}_{l,\epsilon}^c     )  \leq \frac{c}{l^{\frac{\epsilon}{4}}}.
 \end{equation*}

Thus, by the bound \eqref{eqboundness}, we have for all $\epsilon>0$ and $\delta_{l,\epsilon}=\frac{1}{l^{(1-\eta)/2-\epsilon} }$
%
%
%\begin{align*}
%& \P\left(     \sup_{|t-s|\leq \frac{1}{2^l}} |M_n^{\gamma,\beta}[s,t]|  \geq \delta   \right )  \\
%& \leq   \P\left(     \sup_{|u|=l} e^{-\gamma V(u)} ||  M_n(u) ||_{\infty}  \geq \frac{\delta}{3}   \right )  \\   
%& \leq \E \left [     1_{  \sup_{|u|=l} e^{-\gamma V(u)} ||  M_n(u) ||_{\infty}  \geq \frac{\delta}{3} }  | (V(u))_{ |u|=l} \right ]\\ 
%\end{align*}
% 
\begin{align*}
  \P\Big(     \sup_{|t-s|\leq \frac{1}{2^l}} |M_r^{\gamma,\beta}[s,t]|  \geq (\delta_{l,\epsilon})^\gamma  \Big) 
&  \leq   \P( \mathbb{A}_{l,\epsilon}^c     ) +  \P\left(     \sup_{|u|=l} e^{-\gamma V(u)} ||  M^{\gamma,\beta}_r(u) ||_{\infty}  \geq \frac{(\delta_{l,\epsilon})^\gamma}{3} \,  |  \, \mathbb{A}_{l,\epsilon} \right )  \\   
&\leq \frac{c}{l^{\frac{\epsilon}{4}}}+ 1- \P\left(   \max_{|u|=l} \ee^{-\gamma V(u)}   ||  M^{\gamma,\beta}_r(u) ||_{\infty}   \leq  \frac{(\delta_{l,\epsilon})^\gamma}{3} \, | \, \mathbb{A}_{l,\epsilon}  \right)
\\
& \leq  \frac{c}{l^{\frac{\epsilon}{4}}}+ 1- \E\left(  \prod_{|u|=l} \P\left(  \ee^{-\gamma V(u)}    ||  M^{\gamma,\beta}_r(u) ||_{\infty}  \leq  \frac{(\delta_{l,\epsilon})^\gamma}{3}   \:   \big | V(u) \right)     \, | \, \mathbb{A}_{l,\epsilon}        \right).
\end{align*}
In the last line, we have conditioned  on $V(u)$ et used the independence of the mass of the trees rooted at $u$. Hence
\begin{align*}
  \P\Big(     \sup_{|t-s|\leq \frac{1}{2^l}} |M_r^{\gamma,\beta}[s,t]|  \geq \delta_{l,\epsilon}  \Big) 
&    \leq  \frac{c}{l^{\frac{\epsilon}{4}}}+ 1- \E\left(  \prod_{|u|=l}\left\{ 1-  \P\left(     ||  M^{\gamma,\beta}_r(u) ||_{\infty}  \geq  \frac{(\delta_{l,\epsilon})^\gamma}{3} \ee^{\gamma V(u)}  \big | V(u)  \right)       \, | \, \mathbb{A}_{l,\epsilon}    \right\}  \right)
\\
& \leq  \frac{c}{l^{\frac{\epsilon}{4}}}+ 1- \E\left(  \prod_{|u|=l} \left\{1-  \frac{3c}{\delta_{l,\epsilon}}(1+ V(u)^\eta 1_{\min_{j\leq l} V(u_{|j}) \geq -\log l})\ee^{-V(u)}      \right\}      \, | \, \mathbb{A}_{l,\epsilon}     \right)   \\
& \leq  \frac{c}{l^{\frac{\epsilon}{4}}}+   \frac{c}{l^{\frac{\epsilon}{2}}} \\
& \leq  \frac{c}{l^{\frac{\epsilon}{4}}} .
\end{align*}

Therefore, we get
\begin{equation}\label{limsup}
\underset{r \to \infty}{\overline{\lim}}  \P\left(     \sup_{|t-s|\leq \frac{1}{2^l}} |M_r^{\gamma,\beta}[s,t]|  \geq (\delta_{l,\epsilon})^\gamma   \right ) \leq  \frac{c}{l^{\frac{\epsilon}{4}}}
\end{equation}
and then
\begin{equation*}
\underset{l \to \infty}{\lim} \underset{r \to \infty}{\overline{\lim}}  \P\left(     \sup_{|t-s|\leq \frac{1}{2^l}} |M_r^{\gamma,\beta}[s,t]|  \geq (\delta_{l,\epsilon})^\gamma   \right )=0, 
\end{equation*}
which implies that $(M_r^{\gamma,\beta})_{r \geq 1}$ is tight in the space of continuous functions. Now, by using theorem 3 of \cite{kruk} where we view $(M_r^{\gamma,\beta})_{r \geq 1}$ as a martingale which takes values in the Banach space of continuous functions, we conclude that $(M_r^{\gamma,\beta})_{r \geq 1}$ converges almost surely in the space of continuous functions towards some continuous function $M{\gamma,\beta}$. Finally, the estimate \eqref{limsup} implies that
\begin{equation*}
  \P\left(     \sup_{|t-s|\leq \frac{1}{2^l}} |M^{\gamma,\beta}[s,t]|  \geq \delta_{l,\epsilon}   \right ) \leq  \frac{c}{l^{\frac{\epsilon}{4}}}
\end{equation*}
which gives the stated modulus of continuity estimate \eqref{eq:mainmodulus}.
\qed

\section{$M^{\gamma,\beta}_n$ is bounded uniformly in $n$}

The purpose of this section is to prove the main ingredient behind the proof of \ref{th:main}, namely proposition \ref{boundness}. This bound was used precisely in \eqref{thebound} where one can exploit the recursive structure of the tree which entails that $ ||  M^{\gamma,\beta}_r(u) ||_{\infty} $ in  \eqref{thebound} is distributed like the term in the probability on the left hand side of \eqref{eqboundness} (for $p=0$ and $n=r-l$ with $|u|=l$). 

\subsection{Proof of the main  estimate on $M^{\gamma,\beta}$}

We now state the main boundedness estimate of this paper:

\begin{Proposition}
\label{boundness}
There exists $c>0$ and $\eta<1$ such that for any $x\geq 0$, $n\in \N,\, p\geq 0$ 
\begin{equation}\label{eqboundness}
 \P (  ||  M^{\gamma, \beta}_{n,p} ||_{\infty} \geq e^{\gamma x})\leq c(1+x^\eta)\ee^{-x}
\end{equation}
\end{Proposition}
%Proposition \ref{boundness} is proven viaThe following result is a trivial consequence of Proposition \ref{boundness}.

\proof

The proof relies on certain technical lemmas whose proof is postponed to the next sections. We start by introducing for $\epsilon_0\in (0,1)$, $x\geq 0$ the event
\begin{equation}
A(\epsilon_0,x):= \{  \forall u\in \mathbb{T},\ V(u)\geq -x+   r(\epsilon_0) \ln |u|  \}.
\end{equation} 
with $r(\epsilon_0)= \frac{1}{2}-\epsilon_0$. When $x$ is negative and $|x|$ is big, this event will occur with high probability hence one can always work on it: this is the content of lemma \ref{MinAIdimprov}. This is important because we will need tail estimates on $| M^{\gamma,\beta}_{n+p}(z) |$ and in order to get efficient tail estimates, one must condition on $A(\epsilon_0,x)$ (see section \ref{secmoments}). We now fix 
$\epsilon_0$ such that $4\gamma(\frac{1}{2}-\epsilon_0)>1-\epsilon_0$.

 Now, by using lemma \ref{MinAIdimprov}, we get 
 \begin{align*}
& \P\left(  \max_{|z|=n} \left(  \sum_{k=0}^{n-1} \ee^{-\gamma V(z_{|k}^{(l)})} |M_{n+p}^{\gamma,\beta}(z_{|k}^{(l)})|1_{\{ z_{|k+1}\neq z_{|k}^{(l)}  \}} + \ee^{-\gamma V(z)}\left| M^{\gamma,\beta}_{n+p}(z)\right|\right)\geq \ee^{\gamma x}  \right)  \\
& \leq  \P\left(  \max_{|z|=n} \left(  \sum_{k=0}^{n-1} \ee^{-\gamma V(z_{|k}^{(l)})} |M_{n+p}^{\gamma,\beta}(z_{|k}^{(l)})|1_{\{ z_{|k+1}\neq z_{|k}^{(l)}  \}} + \ee^{-\gamma V(z)}\left| M^{\gamma,\beta}_{n+p}(z)\right|\right)\geq \ee^{\gamma x}; \: A(\epsilon_0,x)  \right)  +  c(1+x^{\eta}) \ee^{-x}. \\
\end{align*}
where $\eta \in (0,1)$. Of course, we have
 \begin{align*}
&   \P\left(  \max_{|z|=n} \left(  \sum_{k=0}^{n-1} \ee^{-\gamma V(z_{|k}^{(l)})} |M_{n+p}^{\gamma,\beta}(z_{|k}^{(l)})|1_{\{ z_{|k+1}\neq z_{|k}^{(l)}  \}} + \ee^{-\gamma V(z)}\left| M^{\gamma,\beta}_{n+p}(z)\right|\right)\geq \ee^{\gamma x}; \: A(\epsilon_0,x)  \right)\\
& \leq \P\left(  \max_{|z|=n} \left(  \sum_{k=0}^{n-1} \ee^{-\gamma V(z_{|k}^{(l)})} |M_{n+p}^{\gamma,\beta}(z_{|k}^{(l)})|1_{\{ z_{|k+1}\neq z_{|k}^{(l)}  \}} \right)\geq \ee^{\gamma x}/2; \: A(\epsilon_0,x)  \right)  \\
& + \P\left(  \max_{|z|=n} \left(  \ee^{-\gamma V(z)}\left| M^{\gamma,\beta}_{n+p}(z)\right|\right)\geq \ee^{\gamma x}/2; \: A(\epsilon_0,x)  \right).\\
\end{align*}
Now, we first take care of the term $\P\left(  \max_{|z|=n} \left(  \ee^{-\gamma V(z)}\left| M^{\gamma,\beta}_{n+p}(z)\right|\right)\geq \ee^{\gamma x}/2; \: A(\epsilon_0,x)  \right)$.  
By Lemma \ref{maxondyadique} with $p=2$, $\alpha=r(\epsilon_0)$, $x=x$, $a=x$ and $y=x$, we have

\begin{align*}
& \P\left(  \max_{|z|=n} \left(  \ee^{-\gamma V(z)}\left| M^{\gamma,\beta}_{n+p}(z)\right|\right)\geq \ee^{\gamma x}/2; \: A(\epsilon_0,x)  \right)  \\
& \leq \P\left( \max_{|z|=n} \ee^{-\gamma V(z)}|M_{n+p}^{\gamma,\beta}(z)| 1_{\{   \forall u\geq z,\, V(u)\geq -x+ r(\epsilon_0) \ln |u|  \}} \geq \ee^{\gamma x}  \right) \\
& \leq c(2)n^{r(\epsilon_0)(1-4\gamma)}  \P_x(\underline{S}_n\geq 0) \ee^{-x} \\
& \leq c(2) \ee^{-x}.  \\
\end{align*}

Thus we are left with giving a bound on 
\begin{equation*}
\P\left(  \max_{|z|=n} \left(  \sum_{k=0}^{n-1} \ee^{-\gamma V(z_{|k}^{(l)})} |M_{n+p}^{\gamma,\beta}(z_{|k}^{(l)})|1_{\{ z_{|k+1}\neq z_{|k}^{(l)}  \}} \right)\geq \ee^{\gamma x}/2; \: A(\epsilon_0,x)  \right). 
\end{equation*}
In order to do so, we will split the sum $\sum_{k=0}^{n-1} \ee^{-\gamma V(z_{|k}^{(l)})} |M_{n+p}^{\gamma,\beta}(z_{|k}^{(l)})|1_{\{ z_{|k+1}\neq z_{|k}^{(l)}  \}}$ into two pieces according to the value taken by $V(z_{|k})$. We introduce $\kappa>\frac{40}{\gamma}+1$ and we get  
\begin{align*}
& \P\left(  \max_{|z|=n} \left(  \sum_{k=0}^{n-1} \ee^{-\gamma V(z_{|k}^{(l)})} |M_{n+p}^{\gamma,\beta}(z_{|k}^{(l)})|1_{\{ z_{|k+1}\neq z_{|k}^{(l)}  \}} \right)\geq \ee^{\gamma x}/2; \: A(\epsilon_0,x)  \right) \\
& \leq \P\left(  \max_{|z|=n} \left(  \sum_{k=0}^{n-1} \ee^{-\gamma V(z_{|k}^{(l)})} |M_{n+p}^{\gamma,\beta}(z_{|k}^{(l)})|1_{\{ V(z_{|k}) \leq \kappa \ln k -x \}} \right)\geq \ee^{\gamma x}/2; \: A(\epsilon_0,x)  \right)\\
& + \P\left(  \max_{|z|=n} \left(  \sum_{k=0}^{n-1} \ee^{-\gamma V(z_{|k}^{(l)})} |M_{n+p}^{\gamma,\beta}(z_{|k}^{(l)})|1_{\{ V(z_{|k}) \geq \kappa \ln k -x, \,  z_{|k+1}\neq z_{|k}^{(l)}  \}} \right)\geq \ee^{\gamma x}/2; \: A(\epsilon_0,x)  \right)    \\
\end{align*}

Now, by lemma \ref{PartA}, we get that
\begin{equation*}
 \P\left(  \max_{|z|=n} \left(  \sum_{k=0}^{n-1} \ee^{-\gamma V(z_{|k}^{(l)})} |M_{n+p}^{\gamma,\beta}(z_{|k}^{(l)})|1_{\{ V(z_{|k}) \leq \kappa \ln k -x \}} \right)\geq \ee^{\gamma x}/2; \: A(\epsilon_0,x)  \right) c(1+x^{\eta}) \ee^{-x}
\end{equation*}
The idea behind lemma \ref{PartA} is that typically for very few values of $k$ the event $\{ V(z_{|k}) \leq \kappa \ln k -x \}$ occurs.

Now, by lemma \ref{PartB}, we get that 

\begin{equation*}
\P\left( \max_{|z|=n} \sum_{k=0}^{n-1} \ee^{-\gamma V(z_{|k}^{(l)})} |M_{n+p}^{\gamma,\beta}(z_{|k}^{(l)})|1_{\{    V(z_{|k})\geq \kappa \ln k  -x , \,  z_{|k+1}\neq z_{|k}^{(l)}\}}  \geq \ee^{\gamma x},\, A(x,\epsilon_0)   \right) \leq   c (1+x^\eta)\ee^{-x} .
\end{equation*}

This concludes the proof.

\qed

\subsection{A lemma on the minimum value of all particles}

We first state and prove the following lemma which allows us control the minimum value of all particles:   

\begin{Lemma}
\label{MinAIdimprov} 
Let $\epsilon_0\in (0,\frac{1}{2})$. There exists $c>0$ and $\eta \in (0,1)$ such that for any $x\geq 0$,
\begin{equation}
\P\left( \exists u\in \mathbb{T},\, V(u)\leq -x+ r(\epsilon_0) \ln |u| \right) \leq c(1+x^{\eta}) \ee^{-x}.
\end{equation}
with $r(\epsilon_0) = \frac{1}{2}-\epsilon_0$.
\end{Lemma}
\noindent{\it Proof of Lemma \ref{MinAIdimprov}.} Let $\delta \in (0,1)$. Observe that
\begin{align}
&\nonumber \P\left( \exists u\in \mathbb{T},\ V(u)\leq -x+ r(\epsilon_0) \ln |u|  \right) 
\\
\nonumber &\leq  \E\left( \sum_{j=0}^{+\infty} \sum_{|z|=j}  1_{\{ \min_{i\leq j-1} V(z_{|i})-  r(\epsilon_0)\ln i\geq -x,\, V(z_{|j})-  r(\epsilon_0)\ln j\leq -x  \}}  \right)
\\
\nonumber &\leq  \sum_{j\geq0 }\E\left( \ee^{S_j} 1_{\{  \min_{i\leq j-1} (S_i-  r(\epsilon_0) \ln i) \geq -x,\, S_j\leq -x+  r(\epsilon_0) \ln j  \}}   \right)
\\
\nonumber &\leq \ee^{-x} \sum_{j\geq 0} j^{  r(\epsilon_0)} \P\left( \min_{i\leq j-1} (S_i-  r(\epsilon_0) \ln i) \geq -x,\, S_j\leq -x+  r(\epsilon_0) \ln j   \right)
\\
\nonumber &\leq \ee^{-x} \sum_{j=0}^{x^{2-\delta}} j^{  r(\epsilon_0)} \P\left( S_j\leq -x+  r(\epsilon_0) \ln j   \right)  \\
\nonumber  & +   \ee^{-x} \sum_{j\geq x^{2-\delta} } j^{  r(\epsilon_0)} \P\left( \min_{i\leq j-1} S_i \geq -x,\, S_j\leq -x+  r(\epsilon_0) \ln j   \right)  \\
\nonumber &\leq \ee^{-x} \sum_{j=0}^{x^{2-\delta}} j^{  r(\epsilon_0)} \P\left( S_j\leq -x+  r(\epsilon_0) \ln j   \right)+ c \: \ee^{-x} \sum_{j\geq x^{2-\delta}} j^{  r(\epsilon_0)-\frac{3}{2}} (\ln j)^2(1+x),
\end{align}
where in the last inequality we have used \eqref{estimAid}. By a standard Gaussian estimate, we get
\begin{equation*}
\ee^{-x} \sum_{j=0}^{x^{2-\delta}} j^{  r(\epsilon_0)} \P\left( S_j\leq -x+  r(\epsilon_0) \ln j   \right)  \leq c \ee^{-x-c x^\delta}
\end{equation*}

As $ r(\epsilon_0) -\frac{3}{2}= -1-\epsilon_0<-1$, we deduce that 
\begin{equation*}
\sum_{j\geq x^{2-\delta}} j^{  r(\epsilon_0)-\frac{3}{2}} (\ln j)^2(1+x)\leq c\frac{(1+x)}{ x^{\frac{\epsilon_0}{2}(2-\delta)}}.
\end{equation*}
Gathering the above inequalities leads to
\begin{equation}
\label{ineqAIdBett} \P\left( \exists u\in \mathbb{T},\ V(u)\leq -x+   r(\epsilon_0) \ln |u|  \right)  \leq c(1+x^{\eta}) \ee^{-x},
\end{equation} 
for some $\eta \in (0,1)$. \hfill$\Box$

\subsection{Proof of lemma \ref{PartA}}

Before proving lemma \ref{PartA}, we first prove a large deviation estimate on the number of $k$ such that the event $\{   V(u_{|k})\,  \leq \kappa\ln k -x\} $ occurs. We will see that this number is typically small. Now, for any $x\geq 0$ $\kappa >2$ and $T\in \N^*$, let 
 \begin{equation}
 G(x,\delta, \kappa,T):=  \left\{   \exists  |u|=T^2,\, \exists n_1<...<n_{T^\delta }\in[T,T^2],\, \forall i\in[1,T^{\delta}],\,  V(u_{|n_i})\,  \leq \kappa\ln n_i -x\right\} 
 \end{equation}
 \begin{Lemma}
 \label{MogulModif}
Fix $\kappa>2$. There exists $c_1(\kappa),\, c_2(\kappa)>0$ such that for any $T\geq 10$, $\epsilon_0,\delta >0$, $x\geq 0$,  
\begin{eqnarray}
\label{eqMogulModif}  \P\left( A(x,\epsilon_0)\cap G(x,\delta, \kappa,T)\right)\leq c_1\ee^{- c_2 T^{\frac{\delta}{2}}} \ee^{-x}
\end{eqnarray}
 \end{Lemma}
\noindent{\it Proof of Lemma \ref{MogulModif}.}
The probability in (\ref{eqMogulModif}) is smaller than
\begin{align*}
&  \P\left( A(x,\epsilon_0)\cap G(x,\delta, \kappa,T)\right)
\\
& \leq \E\left(  \sum_{j=T}^{T^2} \sum_{|u|=j}   1_{\{   \inf_{i\leq |u|}(V(u_{|i})-r(\epsilon_0)\ln i)\geq -x, \, V(u)\leq \kappa \ln j -x,\,  \, \exists n_1<...<n_{T^\delta-1}\in [T,j],\, \forall i\in[1,T^{\delta}-1],\,  V(u_{|n_i})   \leq \kappa \ln n_i -x     \}  }      \right)
\\
&\leq \sum_{j=T}^{T^2} \E\left( \ee^{S_j} 1_{\{   \underline{S}_j \geq -x,\,S_j\leq  \kappa \ln j -x, \, \exists n_1<...<n_{T^\delta-1}\in [T,j],\, \forall i\in[1,T^{\delta}],\,  S_{n_i} \leq \kappa \ln n_i -x   \}}   \right)
\\
&\leq \ee^{-x} \sum_{j=T}^{T^2} j^{\kappa} \P\left(    \underline{S}_j \geq -x,\,S_j\leq \kappa \ln j -x,\, \exists n_1<...<n_{T^\delta-1}\in [T,j],\,  \forall i\in[1,T^{\delta}],\,  S_{n_i} \leq \kappa \ln T^2 -x    \right).
\end{align*}
Let 
\begin{eqnarray}
B_j(x,\delta,\kappa) :=  \left\{  \exists n_1<...<n_{T^\delta-1}\in[T,j],\,  \forall i\in[1,T^{\delta}],\,  S_{n_i} \leq \kappa \ln (T^2)  -x  \right\},
\end{eqnarray}
we claim that there exists $c_1,c_2>0$ such that for any $j\in [T,T^2]$, 
\begin{eqnarray}
\label{clami}
\P\left(\{ \underline{S}_j\geq -x,\, S_j\leq \kappa \ln (T^2)-x \} \cap  B_j(x,\delta ,\kappa)  \right) \leq c_1\ee^{-c_2 T^{\frac{\delta}{2}}}.
\end{eqnarray}
Inequality (\ref{clami}) is sufficient to prove Lemma \ref{MogulModif}, indeed assuming (\ref{clami}) we get
\begin{eqnarray}
\P\left( A(x,\epsilon_0)\cap G(x,\delta, \kappa,T)\right)   \leq \sum_{j=T}^{T^2} \ee^{-x} T^{2\kappa} c_1\ee^{-c_2 T^{\delta\over 2}} \leq \ee^{-x} c_1'\ee^{-c_2' T^{\delta\over 2}}.
\end{eqnarray}
Thus it remains to prove \eqref{clami}.

\begin{figure}[h]
\centering
\caption{}
\label{intttux2}
\includegraphics[interpolate=true,width=15cm,height=1.5cm]{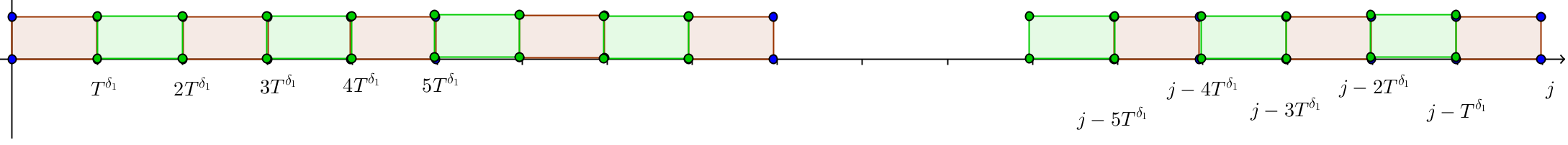} 
\end{figure}

Let $\delta_1\in (0,\delta)$, say $\delta_1={\delta\over 2}$. There exists at least $T^{\delta-\delta_1}= T^{\delta_1}$ red or green (see figure 1) intervals which contain one $n_i$. For instance we can assume that the red intervals contain at least $t(\delta):=\frac{T^{\delta\over 2}}{2}$ times  $n_i$. We define the sequence of stopping times:
\begin{eqnarray*}
\sigma_0&:=&0,
\\
\sigma_1&:=& \inf\{i> t(\delta),\, S_i\leq \kappa 2\ln T-x   \},
\\
\tau_1 &:=&\inf \{  i\geq \sigma_1+ t(\delta),\, \exists k\in \N,\, i=kt(\delta)\}
\\
\sigma_2 &:=& \inf\{ i\geq \tau_1,\, S_i\leq \kappa2\ln T-x   \}
\\
\tau_2 &:= & \inf\{  i\geq \sigma_2+ t(\delta),\, \exists k\in \N,\, i=kt(\delta)  \},
\end{eqnarray*}
and so on. On the event $  B_j(x,\delta ,\kappa) $ we have $\sigma_{t(\delta)-2}<j-T^{\delta\over 2}$. So we shall bound 
\begin{equation}
\P\left( \underline{S}_j\geq -x,\, \sigma_{t(\delta)-1}<j- T^{\delta\over 2},\, S_j\leq \kappa 2\ln T-x    \right)
\end{equation}
Moreover by the Markov property at time $\sigma_{t(\delta)-1}$ we get
\begin{eqnarray*}
&&\P\left( \underline{S}_j\geq -x,\, \sigma_{t(\delta)-1}<j- T^{\delta\over 2},\, S_j\leq \kappa 2\ln T-x    \right)
\\
&&\leq \E\left( 1_{\{ \underline{S}_{\sigma_{t(\delta)-1}}\geq -x,\, \sigma_{t(\delta)-1}\leq j- T^{\delta\over 2} \}} \P_{S_{\sigma_{t(\delta)-1}} + x}\left(  \underline{S}_{j-\sigma_{t(\delta)-1}}\geq 0,\, S_{j-\sigma_{t(\delta)-1}}\leq 2\kappa\ln T      \right)   \right)
\\
&&\leq \frac{c (2\kappa \ln T)^3}{(T^{\delta\over 2})^{\frac{3}{2}}}\P\left(  \underline{S}_{\sigma_{t(\delta)-1}}\geq -x,\, \sigma_{t(\delta)-1}\leq j- T^{\delta\over 2}  \right),
\end{eqnarray*}
where in the last line we used the fact that $j-\sigma_{t(\delta)-1}>T^{\delta\over 2}$. Now by the Markov property at time $\sigma_{t(\delta)-2}$ we get
\begin{eqnarray*}
&&\P\left( \underline{S}_j\geq -x,\, \sigma_{t(\delta)-1}<j-  T^{\delta\over 2},\, S_j\leq \kappa 2\ln T-x    \right)
\\
&&\leq  \frac{c (2\kappa \ln T)^3}{(  T^{\delta\over 2}  )^{\frac{3}{2}}} \E\left(  1_{\{ \underline{S}_{   \sigma_{t(\delta)-2}    }\geq -x,\, \sigma_{  t(\delta)-2   }\leq j-2  T^{\delta\over 2} \}} \times \sup_{y\in [-x,2\kappa\ln T-x]} \sum_{i\geq  T^{\delta\over 2}} \P_y\left(  \underline{S}_i\geq -x,\, S_i\leq 2\kappa\ln T-x \right)     \right)
\\
&&\leq   \frac{c (2\kappa \ln T)^3}{(  T^{\delta\over 2})^{\frac{3}{2}}} \frac{c (\kappa2\ln T)^3}{(  T^{\delta\over 2})^\frac{1}{2}} \P\left(  \underline{S}_{   \sigma_{t(\delta)-2}    }\geq -x,\, \sigma_{  t(\delta)-2   }\leq j- 2  T^{\delta\over 2}   \right)
\end{eqnarray*}
By iterating yet $t(\delta)-2$ times this procedure we get
\begin{equation}
\P(  \{ \underline{S}_j\geq -x,\, S_j\leq 2\ln (T^2)-x \} \cap  B_j(x,\delta ,\kappa) )\leq  \frac{c (2\kappa \ln T)^3}{(  T^{\delta\over 2})^{\frac{3}{2}}}  \left( \frac{c (\kappa2\ln T)^3}{(t(\delta))^\frac{1}{2}}\right)^{t(\delta)-1}\leq c_1\exp(-c_2  T^{\delta\over 2}) ,
\end{equation} 
which concludes the proof of the Lemma \ref{MogulModif}. \hfill$\Box$

Now, we can state the main lemma of this subsection:

\begin{Lemma}
\label{PartA}
For any $\kappa>0$ and $\epsilon_0>0$ such that $4\gamma(\frac{1}{2}-\epsilon_0)>1-\epsilon_0$, there exists $c>0$ and $\eta \in (0,1)$ such that for any $x\geq 0,\, n\in \N$, 
\begin{equation*}
\P\left( \max_{|z|=n} \sum_{k=0}^{n-1} \ee^{-\gamma V(z_{|k}^{(l)})} |M_{n+p}^{\gamma,\beta}(z_{|k}^{(l)})|1_{\{  V(z_{|k})\leq \kappa \ln k  -x \}}  \geq \ee^{\gamma x},\, A(x,\epsilon_0) \right) \leq c (1+x^\eta) \ee^{-x} .
\end{equation*}
\end{Lemma}
\noindent{\it Proof of Lemma \ref{PartA}.}  First observe that on $A(x,\epsilon_0)$ the condition $  \min_{j\leq k+1} {V}((z ^{(l)}) _{|j}) \geq -x$ is automatically satisfied for any $z\in \mathbb{T}$. For any $x,\, y,\, a\geq 0$ and $r\in (0,\frac{1}{2})$, we define
\begin{equation}\label{defdeO}
O_k(x,y,a,r):= \{       \max_{|z|=k} \ee^{-\gamma V(z)} 1_{\{ \min_{j\leq k}V(z_{|j})\geq -a\}} |M_{k+n}^{\gamma,\beta}(z)| 1_{\{   \forall u\geq z,\, V(u)\geq -x+ r \ln |u|  \}} \geq \ee^{\gamma y}        \}
\end{equation}
Since $4\gamma(\frac{1}{2}-\epsilon_0-\theta)>1-\epsilon_0$, one can choose $\theta>0$ small enough such that $ 4\gamma r(\epsilon_0) - 4\gamma \theta- r(\epsilon_0)>\frac{1}{2}$. By applying Lemma \ref{maxondyadique}, with $x=x$, $y=x-\theta\ln k$, $a=x$, $r=r(\epsilon_0)$, we get
\begin{eqnarray}
\nonumber \sum_{k\geq 1} \P\left( O_k(x,x-\theta\ln k, x, r(\epsilon_0))  \right)&\leq& \sum_{k\geq 0}   c(2) k^{r(\epsilon_0)(1-4 \gamma)} \ee^{4\gamma \theta \ln k} \ee^{-x}\P_x(\underline{S}_k\geq 0)
\\
\nonumber &\leq &  c \ee^{-x}  \sum_{k\geq 1}  k^{r(\epsilon_0)+ 4\gamma \theta -4\gamma r(\epsilon_0) }   \P_x(\underline{S}_k\geq 0)
\\
\nonumber &\leq &  c \ee^{-x}  \sum_{k\geq 1}  k^{r(\epsilon_0)+ 4\gamma \theta -4\gamma r(\epsilon_0) }   \P(|N| \leq \frac{x}{\sqrt{k}})
%\\
%&\leq & c(1+x^\eta) \ee^{-x}.
\end{eqnarray}
where $N$ is a standard Gaussian variable. Now, we have that
\begin{align}
 & \nonumber  \ee^{-x}  \sum_{k\geq 1}  k^{r(\epsilon_0)+ 4\gamma \theta -4\gamma r(\epsilon_0) }   \P(|N| \leq \frac{x}{\sqrt{k}})
\\
& \nonumber \leq   c \ee^{-x}  \sum_{k=1}^x k^{r(\epsilon_0)+ 4\gamma \theta -4\gamma r(\epsilon_0) }   \P(|N| \leq \frac{x}{\sqrt{k}}) + c \ee^{-x}  \sum_{k \geq x} k^{r(\epsilon_0)+ 4\gamma \theta -4\gamma r(\epsilon_0) }   \P(|N| \leq \frac{x}{\sqrt{k}})  \\
&\nonumber  \leq   c \ee^{-x}  \sqrt{x} + c \ee^{-x} (1+x)  \sum_{k \geq x} k^{r(\epsilon_0)+ 4\gamma \theta -4\gamma r(\epsilon_0) -\frac{1}{2} }  
\\
& \nonumber  \leq  c(1+x^\eta) \ee^{-x}
\end{align}
for some $\eta\in (0,1)$. Hence, in conclusion, we get that
\begin{equation*}
\sum_{k\geq 1} \P\left( O_k(x,x-\theta\ln k, x, r(\epsilon_0))  \right)  \leq c(1+x^\eta) \ee^{-x}
\end{equation*}

%Thus to prove Lemma \ref{PartA}, it suffices to prove that
%\begin{eqnarray*}
%\P\left(  \max_{|z|=n} \sum_{k=0}^{n-1} \frac{1}{k^{\theta \gamma}} 1_{\{ \underline{V}(z_k^{(l)})\geq -x,\,   V(z_k)\leq \kappa\ln k  -x \}} \geq 1 ,\, A(x,\epsilon_0)   \right)\leq c\ee^{-x}.
%\end{eqnarray*}
Observe that for any $|z|=n$, on the set  $\cap_{p\geq 0} \left(  O_k(x,x-\theta\ln k, x, r(\epsilon_0))^c  \right) \cap A(x,\epsilon_0) $, we have
\begin{align*}
 & \sum_{k=T}^{n-1} \ee^{-\gamma V(z_{|k}^{(l)})} |M_{n+p}^{\gamma,\beta}(z_{|k}^{(l)})|1_{\{  V(z_{|k})\leq \kappa \ln k  -x \}}  \\
  & = \sum_{k=T}^{n-1} \ee^{-\gamma V(z_{|k}^{(l)})} |M_{n+p}^{\gamma,\beta}(z_{|k}^{(l)})|1_{\{  \min_{j\leq k+1} {V}((z_{|k}^{(l)}) _{|j})\geq -x,\, \forall u\geq z_{|k}^{(l)},\, V(u)\geq r(\epsilon_0)\ln |u|-x,\, V(z_{|k})\leq \kappa \ln k -x  \}}  \\
  &\leq   \ee^{\gamma x} \sum_{k=T}^{n-1} \frac{1}{k^{\theta \gamma}} 1_{\{   \min_{j\leq k+1} {V}((z_{|k}^{(l)}) _{|j})\geq -x,\,   V(z_{|k})\leq \kappa\ln k  -x \}}
\end{align*}
Let $\delta = \frac{\theta\gamma}{2}$. By lemma \ref{MogulModif}, recall that 
\begin{equation}
\sum_{p\geq 0} \P\left( A(x,\epsilon_0)\cap   G(x,\delta,\kappa, T^{2^{p}})  \right)\leq \sum_{p\geq 0} \ee^{-x}c_1 \ee^{-c_2 T^{2^{p}\delta\over 2}} \leq  c\ee^{-x}
\end{equation}
Moreover on the set $\cap_{p\geq 0} \left( G(x,\delta,\kappa,T^{2^p})^c   \right)$, we have
\begin{align*}
&\sup_{|z|=n}\sum_{k=T}^{n-1} \frac{1}{k^{\theta \gamma}} 1_{\{  \min_{j\leq k+1} {V}((z_{|k}^{(l)}) _{|j})  \geq -x,\,   V(z_{|k})\leq \kappa\ln k  -x     \}} \\
&   \leq  \sup_{|z|=n}   \sum_{p\geq 0} \sum_{k=T^{2^{p}}}^{\min(T^{2^{p+1}}, n)} \frac{1}{k^{\theta\gamma}} 1_{\{ \min_{j\leq k+1} {V}((z_{|k}^{(l)}) _{|j})  \geq -x,\,   V(z_{|k})\leq \kappa\ln k  -x \}}   \\
& \leq   \sum_{p\geq 0}\frac{T^{2^p\delta} }{T^{2^p\theta\gamma}}  \\
& \leq c(\kappa, \theta\gamma,T) 
\end{align*}
We choose  $T$ large enough such that $c(\kappa, \theta\gamma,T)<1$ and we set $$\mathcal{B}=\cap_{p\geq 0} \left(  O_k(x,x-\theta\ln k, x, r(\epsilon_0))^c  \right) \cap \cap_{p\geq 0} \left( G(x,\delta,\kappa,T^{2^p})^c   \right) \cap A(x,\epsilon_0) . $$ Therefore, on the event $\mathcal{B}$, we have
\begin{equation*}
\max_{|z|=n} \sum_{k=0}^{n-1} \ee^{-\gamma V(z_{|k}^{(l)})} |M_{n+p}^{\gamma,\beta}(z_{|k}^{(l)})|1_{\{  V(z_{|k})\leq \kappa \ln k -x  \}} \leq \max_{|z|=T} \sum_{k=0}^{T} \ee^{-\gamma V(z_{|k}^{(l)})} |M_{n+p}^{\gamma,\beta}(z_{|k}^{(l)})|1_{\{  V(z_{|k})\leq \kappa \ln k -x  \}} + \ee^{\gamma x} c(\kappa, \theta\gamma,T)
\end{equation*}

We deduce that
\begin{align*}
&  \P\left( \max_{|z|=n} \sum_{k=0}^{n-1} \ee^{-\gamma V(z_{|k}^{(l)})} |M_{n+p}^{\gamma,\beta}(z_{|k}^{(l)})|1_{\{ \forall u\geq z_{|k}^{(l)},\, V(u)\geq r(\epsilon_0)\ln |u|-x,\, V(z_{|k})\leq \kappa \ln k -x  \}}  \geq \ee^{\gamma x},\, A(x,\epsilon_0) \right)
\\
&\leq c(1+x^\eta)\ee^{-x} +\P\left( \max_{|z|=T} \sum_{k=0}^{T} \ee^{-\gamma V(z_{|k}^{(l)})} |M_{n+p}^{\gamma,\beta}(z_{|k}^{(l)})|  \geq \ee^{\gamma x}[1 - c(\kappa, \theta\gamma,T)],\, \mathcal{B}  \right)
\\
&\leq c(1+x^\eta) \ee^{-x}.
\end{align*}
where in the last line we have used Lemma \ref{momentp2}. 

\subsection{Proof of lemma \ref{PartB}}

\begin{Lemma}
\label{PartB}
Let $\kappa>  \frac{40}{\gamma}+1$  and $\epsilon_0$ be such that $4\gamma(\frac{1}{2}-\epsilon_0)>1-\epsilon_0$. There exists $c>0$ such that for any $x \geq 0,\, n,\, p\in \N$, 
\begin{eqnarray*}
&&\P\left( \max_{|z|=n} \sum_{k=0}^{n-1} \ee^{-\gamma V(z_{|k}^{(l)})} |M_{n+p}^{\gamma,\beta}(z_{|k}^{(l)})|1_{\{    V(z_{|k})\geq \kappa \ln k  -x, \: z_{|k+1}\neq z_{|k}^{(l)} \}}  \geq \ee^{\gamma x},\, A(x,\epsilon_0)   \right) \leq   c (1+x^\eta)\ee^{-x} .
\end{eqnarray*}

\end{Lemma}
\noindent{\it Proof of Lemma \ref{PartB}.} Let $\upsilon\in(0,1)$, such that $2- (\kappa-2-r(\epsilon_0)) \frac{\gamma\upsilon}{2}<0$ and $4\gamma(1-\upsilon) -1>0$. Such a $\upsilon$ exists since $\kappa >\frac{40}{\gamma}+1 $. Let $\epsilon_0$ such that $4\gamma r(\epsilon_0)-r(\epsilon_0)>\frac{1}{2}$.   Recall definition \eqref{defdeO} which was introduced in the proof of lemma \ref{PartA} 
\begin{equation*}
O_k(x,y,a,r)= \{       \max_{|z|=k} \ee^{-\gamma V(z)} 1_{\{ \min_{j\leq k}V(z_{|j})\geq -a\}} |M_{k+n}^{\gamma,\beta}(z)| 1_{\{   \forall u\geq z,\, V(u)\geq -x+ r \ln |u|  \}} \geq \ee^{\gamma y}        \}
\end{equation*}

We have
\begin{eqnarray}
\nonumber \sum_{t=1}^{+\infty }\sum_{k=1}^{+\infty} \P(   O_k(x-t,x-\upsilon t,x, r(\epsilon_0)),\, A(x,\epsilon_0)) &\leq& \sum_{t=1}^{\infty}\sum_{k=1}^{+\infty} c(2) k^{r(\epsilon_0)(1-4 \gamma)} \ee^{4\gamma (\upsilon t-t)} \ee^{-x}\P_x(\underline{S}_k\geq 0) \ee^{t}
\\
\nonumber &\leq &c(1+x^{\eta})\ee^{-x}\sum_{t=1}^{\infty} \ee^{-4\gamma (1-\upsilon)t +t}
\\
\label{touse} &\leq & c(1+x^{\eta})\ee^{-x},
\end{eqnarray}
where we have used the fact (which appears in the proof of lemma \ref{PartA}) that $$\sum_{k=1}^{+\infty} k^{r(\epsilon_0)(1-4 \gamma)} \P_x(\min_{j\leq k}{S}_j\geq 0) \leq c(1+x^{\eta}). $$
Otherwise it is clear that on $A(x,\epsilon_0)$,
\begin{align*}
 &\max_{|z|=n} \sum_{k=0}^{n-1} \ee^{-\gamma V(z_{|k}^{(l)})} |M_{n+p}^{\gamma,\beta}(z_{|k}^{(l)})|1_{\{ V(z_{|k})\geq \kappa \ln k -x , \: z_{|k+1}\neq z_{|k}^{(l)}  \}}  
\\
& = \max_{|z|=n} \sum_{k=0}^{n-1} \ee^{-\gamma V(z_{|k}^{(l)})} |M_{n+p}^{\gamma,\beta}(z_{|k}^{(l)})|1_{\{  \min_{j\leq k+1} {V}((z_{|k}^{(l)}) _{|j})\geq -x,\, \forall u\geq z_{|k}^{(l)},\, V(u)\geq r(\epsilon_0)\ln |u|-x,\, V(z_{|k})\geq \kappa \ln k -x , \: z_{|k+1}\neq z_{|k}^{(l)}  \}}   
\\
& \leq \max_{|z|=n} \sum_{k=0}^{n-1} \sum_{t=0}^{+\infty}\frac{|M_{n+p}^{\gamma,\beta}(z_{|k}^{(l)})|}{ \ee^{\gamma V(z_{|k}^{(l)})}} \\  
& \times 1_{\{V(z_{|k})\geq \kappa\ln k-x,\,   \min_{j\leq k+1} {V}((z_{|k}^{(l)}) _{|j})\geq -x,\,  \forall u\geq z_{|k}^{(l)},\, V(u)\geq r(\epsilon_0)\ln |u|-x+t,\, \exists u\geq z_{|k}^{(l)},\, V(u)<r(\epsilon_0)\ln |u|-x+t+1 , \: z_{|k+1}\neq z_{|k}^{(l)} \}}.
\end{align*}
Let us set $\bar{\mathcal{B}}= \cap_{t \geq 0} \cap_{k \geq 0} O_k(x-t,x-\upsilon t,x, r(\epsilon_0)) ^c$. By using (\ref{touse}) and the definition of $O_k$ , we deduce that 
\begin{align*}
&\P\left( \max_{|z|=n} \sum_{k=0}^{n-1} \ee^{-\gamma V(z_{|k}^{(l)})} |M_{n+p}^{\gamma,\beta}(z_{|k}^{(l)})|1_{\{   \forall u\geq z_{|k}^{(l)},\, V(u)\geq r(\epsilon_0)\ln |u|-x,\, V(z_{|k})\geq \kappa \ln k -x , \: z_{|k+1}\neq z_{|k}^{(l)} \}}  \geq \ee^{\gamma x} ,\, A(x,\epsilon_0) \right)
\\
&\leq c(1+x^{\eta})\ee^{-x}+    \P\left(    \max_{|z|=n} \sum_{k=0}^{n-1} \sum_{t=0}^{+\infty} \ee^{-\gamma \upsilon t} 1_{\{  \exists u\geq z_{|k}^{(l)},\, V(u)<r(\epsilon_0)\ln |u|-x+t+1 ,  \, V(z_{|k})\geq \kappa\ln k-x , \, z_{|k+1}\neq z_{|k}^{(l)} \}}        \geq 1   ,\, A(x,\epsilon_0) \cap \bar{\mathcal{B}}   \right)
\\
&\leq c(1+x^{\eta})\ee^{-x}+    \P\left(  \sum_{t=0}^{+\infty}  \max_{|z|=n} \sum_{k=0}^{n-1}  \ee^{-\gamma \upsilon t} 1_{\{  \exists u\geq z_{|k}^{(l)},\, V(u)<r(\epsilon_0)\ln |u|-x+t+1 ,  \, V(z_{|k})\geq \kappa\ln k-x , \, z_{|k+1}\neq z_{|k}^{(l)} \}}        \geq 1   ,\, A(x,\epsilon_0) \cap \bar{\mathcal{B}}   \right)
\\
&\leq  c(1+x^{\eta})\ee^{-x}  \\
& + \sum_{t=0}^{+\infty} \P\left(     \max_{|z|=n} \sum_{k=0}^{n-1}  1_{\{  \exists u\geq z_{|k}^{(l)},\, V(u)<r(\epsilon_0)\ln |u|-x+t+1, \, V(z_{|k})\geq \kappa\ln k-x ,\, \, z_{|k+1}\neq z_{|k}^{(l)} \}}        \geq   (1-\ee^{\gamma\upsilon\over 2})\ee^{\gamma \upsilon t\over 2}     ,\, A(x,\epsilon_0)    \right)
\end{align*}
We set $c_{\upsilon \gamma}:= \frac{1}{2}(1-\ee^{\gamma\upsilon\over 2})$. We deduce that 
\begin{align*}
&\P\left( \max_{|z|=n} \sum_{k=0}^{n-1} \ee^{-\gamma V(z_{|k}^{(l)})} |M_{n+p}^{\gamma,\beta}(z_{|k}^{(l)})|1_{\{   \forall u\geq z_{|k}^{(l)},\, V(u)\geq r(\epsilon_0)\ln |u|-x,\, V(z_{|k})\geq \kappa \ln k -x , \: z_{|k+1}\neq z_{|k}^{(l)} \}}  \geq \ee^{\gamma x} ,\, A(x,\epsilon_0) \right)	
\\
&\leq c (1+x^{\eta}) \ee^{-x}  \\
&+   \sum_{t=  0}^{+\infty} \P\left(     \max_{|z|=n} \sum_{k=  c_{\upsilon \gamma} \ee^{\gamma \upsilon t\over 2} }^{n-1} 1_{\{  \exists u\geq z_{|k}^{(l)},\, V(u)<r(\epsilon_0)\ln |u|-x+t+1 ,\, V(z_{|k})\geq \kappa \ln k -x , \: z_{|k+1}\neq z_{|k}^{(l)} \}}        \geq  c_{\upsilon \gamma} \ee^{\gamma \upsilon t\over 2}       ,\, A(x,\epsilon_0)  \right). 
\end{align*}

Let $|z|=n$ be such that $ \sum_{k=   c_{\gamma\upsilon}  \ee^{\gamma \upsilon t\over 2} }^{n-1}   1_{\{   \exists u\geq z_{|k}^{(l)},\, V(u)<r(\epsilon_0)\ln |u|-x+t+1, \, V(z_{|k})\geq \kappa\ln k-x , \: z_{|k+1}\neq z_{|k}^{(l)}  \}} \neq 0$. If $\tilde{k}$ is the last index such that $1_{\{ \exists u\geq z_{| \tilde{k}}^{(l)},\, V(u)<r(\epsilon_0)\ln |u|-x+t+1, \, V(z_{\tilde{k}})\geq \kappa\ln {\tilde{k}}-x,\, \: z_{|\tilde{k}+1}\neq z_{|\tilde{k}}^{(l)} \}}\neq 0$ then there exists $v \geq z_{|\tilde{k}}^{(l)}$ such that 
\begin{equation}
V(v) <r(\epsilon_0) \ln |v|-x+t+1.
\end{equation}
Notice that $p=|v|\geq \tilde{k} \geq c_{\gamma\upsilon}  \ee^{\gamma \upsilon t\over 2}$ and that, by definition of $\tilde{k}$, we have
\begin{align*}
& \sum_{k=  c_{\upsilon \gamma} \ee^{\gamma \upsilon t\over 2} }^{n-1} 1_{\{  \exists u\geq z_{|k}^{(l)},\, V(u)<r(\epsilon_0)\ln |u|-x+t+1 ,\, V(z_{|k})\geq \kappa \ln k -x , \: z_{|k+1}\neq z_{|k}^{(l)} \}}          \\
& \leq     \sum_{k=  c_{\upsilon \gamma} \ee^{\gamma \upsilon t\over 2} }^{|v|} 1_{\{  \exists u\geq v_{|k}^{(l)},\, V(u)<r(\epsilon_0)\ln |u|-x+t+1 ,\, V(v_{|k})\geq \kappa \ln k -x , \: v_{|k+1}\neq v_{|k}^{(l)} \}}    \\
\end{align*}
Therefore we get for any $|z|=n$, 
\begin{eqnarray*}
&&\sum_{k=   c_{\upsilon \gamma}  \ee^{ \upsilon \gamma t\over 2} }^{n-1}   1_{\{   \exists u\geq z_{|k}^{(l)},\, V(u)<r(\epsilon_0)\ln |u|-x+t+1 ,\, V(z_{|k})\geq \kappa \ln k -x , \: z_{|k+1}\neq z_{|k}^{(l)} \}} 
 \\
&& \leq \sum_{p=  c_{\upsilon \gamma}  \ee^{ \upsilon \gamma t\over 2} }^{+\infty}    \sum_{|z|=p}  1_{\{    V(z)\leq r(\epsilon_0)\ln p+ t-x  \}}  \left(  \sum_{k=  c_{\upsilon \gamma}  \ee^{\gamma \upsilon t\over 2}}^{\min(p,n-1)   }   1_{\{    \exists u\geq z_{|k}^{(l)},\, V(u)<r(\epsilon_0)\ln |u|-x+t+1 ,\, V(z_{|k})\geq \kappa \ln k -x , \: z_{|k+1}\neq z_{|k}^{(l)} \}} \right)   
\\
&&
\leq   \sum_{p=  c_{ \upsilon \gamma} \ee^{\frac{\upsilon \gamma t}{2}} }^{+\infty}    \sum_{|z|=p}  1_{\{    V(z)\leq r(\epsilon_0)\ln p+ t-x  \}}  \left(  \sum_{k=  c_{\upsilon \gamma}  \ee^{\gamma \upsilon t\over 2}}^{\min(p,n-1)   }   1_{\{    \exists u\geq z_{|k}^{(l)},\, V(u)<r(\epsilon_0)\ln |u|-x+t+1 ,\, V(z_{|k})\geq \kappa \ln k -x , \: z_{|k+1}\neq z_{|k}^{(l)} \}} \right)
\end{eqnarray*}
Moreover on $A(x,\epsilon_0)$, we can add the indicator function $ 1_{\{ \min_{j\leq |z|}(V(z_{|j})-r(\epsilon_0)\ln j)\geq  -x  \}} $ which leads to 
\begin{align*}
& \P\left( \max_{|z|=n} \sum_{k=   c_{\upsilon \gamma}  \ee^{\gamma \upsilon t\over 2} }^{n-1}   1_{\{    \exists u\geq z_{|k}^{(l)},\, V(u)<r(\epsilon_0)\ln |u|-x+t+1 ,\, V(z_{|k})\geq \kappa \ln k -x , \: z_{|k+1}\neq z_{|k}^{(l)} \}}  \geq   c_{\upsilon \gamma} \ee^{\gamma \upsilon t\over 2}    ,\, A(x,\epsilon_0) \right)
\\
& \leq  \P\left(  \sum_{p=  c_{ \upsilon \gamma} \ee^{\frac{\upsilon \gamma t}{2}} }^{+\infty}    \sum_{|z|=p}  1_{\{    V(z)\leq r(\epsilon_0)\ln p+ t-x  ,\,  \min_{j\leq |z|}(V(z_{|j})-r(\epsilon_0)\ln j)\geq  -x \}} \times \right.
\\
& \qquad\qquad\qquad \qquad \left. \left(  \sum_{k=  c_{\upsilon \gamma}  \ee^{\gamma \upsilon t\over 2}}^{\min(p,n-1)   }   1_{\{   \exists u\geq z_{|k}^{(l)},\, V(u)<r(\epsilon_0)\ln |u|-x+t+1 ,\, V(z_{|k})\geq \kappa \ln k -x , \: z_{|k+1}\neq z_{|k}^{(l)} \}} \right)  \geq   1  \right)
\\
&\leq  \sum_{p= c_{\upsilon \gamma}  \ee^{\gamma \upsilon t\over 2} }^{+\infty}   \sum_{|z|=p} \E\left(  1_{\{ \min_{j\leq |z|}(V(z_{|j})-r(\epsilon_0)\ln j)\geq  -x,\,  V(z)\leq r(\epsilon_0)\ln p+ t-x  \}} \times\right. 
 \\
 & \qquad\qquad\qquad\qquad\qquad  \left. \left(  \sum_{k=  c_{\upsilon \gamma}  \ee^{\gamma \upsilon t\over 2}}^{\min(p,n-1)   }   1_{\{   \exists u\geq z_{|k}^{(l)},\, V(u)<r(\epsilon_0)\ln |u|-x+t+1 ,\, V(z_{|k})\geq \kappa \ln k -x , \: z_{|k+1}\neq z_{|k}^{(l)} \}} \right)  \right)
 \end{align*}

% \\
% &&\leq \sum_{p= c_{\gamma\eta} \ee^{\gamma \eta t\over 2} }^{+\infty} (1+x^{\eta})\ee^{-x+t}\frac{p^{r(\epsilon_0)}}{p^{\frac{3}{2}-\frac{\epsilon_0}{2}}} \E_\Q\left( 1_{\{ \min_{j\leq |z|}(V(z_j)-r(\epsilon_0)\ln j)\geq  -x,\,  V(z)\leq r(\epsilon_0)\ln p+ t-x  \}}  \sum_{j= c_{\gamma\eta} \ee^{\gamma \eta t\over 2}}^{\min(p,n-1)   }   1_{\{  V(w_j)\geq \kappa\ln j-x,\,  \exists u\geq w_k^{(l)},\, V(u)<r(\epsilon_0)\ln |u|-x+t+1 \}}   \right),  
% \\
%&&\leq \sum_{p= c_{\gamma\eta} \ee^{\gamma \eta t\over 2} }^{+\infty} \ee^{-x+t}\frac{p^{r(\epsilon_0)}}{p^{\frac{3}{2}-\frac{\epsilon_0}{2}}}  \E\left(  \sum_{|z|=p}  \left(  \sum_{j=  c_{\gamma\upsilon}  \ee^{\gamma \upsilon t\over 2}}^{\min(p,n-1)   }   1_{\{  V(z_j)\geq \kappa\ln j-x,\,  \exists u\geq z_k^{(l)},\, V(u)<r(\epsilon_0)\ln |u|-x+t+1 \}} \right)\right ),
%\end{eqnarray*}
%with $w$ the trajectory associated to the spine under $\Q$. 

From Lemma \ref{MinAIdimprov}, recall that for any $y\geq 0$
\begin{equation}
\label{recall} \P\left( \exists v\in \mathbb{T},\, V(v)\leq r(\epsilon_0)\ln |v|-y\right)\leq c(1+y^{\eta})\ee^{-y}.
\end{equation}
By taking the conditional expectation according to the sigma field $\sigma( V(z_{|k}),\, k\leq p):= \mathcal{G}_z$, via the branching property, by using the inequality $\ln |u|\leq \ln k+ \ln (|u|-k)$ for any $|u|\geq j$, and applying (\ref{recall}) with $y=\kappa \ln j-r(\epsilon_0)\ln j-t-1$,  we get 
\begin{align*}
&\E\left(  1_{\{   \exists u\geq z_{|k}^{(l)},\, V(u)<r(\epsilon_0)\ln |u|-x+t+1 ,\, V(z_{|k})\geq \kappa \ln k -x , \: z_{|k+1}\neq z_{|k}^{(l)} \}} \big| \mathcal{G}_z \right)\leq \E\left( 1_{\{V(z_{|k}^{(l)})-V(z_{|k}) \geq \ln k   \}} \right)
\\
&     \E\left(   \left[  c(1+ (\kappa\ln k-t- r(\epsilon_0) \ln k- [V(z_{|k}^{(l)})-V(z_{|k})]  )^{\eta}  ) \ee^{-  (  \kappa\ln k-t- r(\epsilon_0) \ln k-[V(z_{|k}^{(l)})-V(z_{|k})] )}    \right] 1_{\{V(z_{|k}^{(l)})-V(z_{|k}) \leq \ln k    \}}   \right)
\\
&\leq  c\ee^{-(\ln k)^2}  +  \left[  c(1+ ((2 \times \frac{2}{\gamma \upsilon}+1) \ln k -t )^{\eta}  ) \ee^{-  ( (2 \times \frac{2}{\gamma \upsilon}+1) \ln k -t )}    \right]   
\end{align*}
where we have used $(\kappa-2-r(\epsilon_0))\frac{\gamma n}{2}>2$ in the last line. Finally we deduce that 
\begin{align*}
	&\P\left( \max_{|z|=n} \sum_{k=   c_{\gamma\upsilon}  \ee^{\gamma \upsilon t\over 2} }^{n-1}   1_{\{  V(z_{|k})\geq \kappa\ln k-x,\,  \exists u\geq z_{|k}^{(l)},\, V(u)<r(\epsilon_0)\ln |u|-x+t+1  \}}  \geq   c_{\gamma\upsilon} \ee^{\gamma \upsilon t\over 2}    ,\, A(x,\epsilon_0) \right)
 \\
& \sum_{p= c_{\gamma\upsilon}  \ee^{\gamma \upsilon t\over 2} }^{+\infty} \E\left(  \sum_{|z|=p}   1_{\{ \min_{j\leq |z|}(V(z_{|j})-r(\epsilon_0)\ln j)\geq  -x,\,  V(z)\leq r(\epsilon_0)\ln p+ t-x  \}}  \right)\times  
  \\
  &   \qquad\qquad \qquad    \sum_{k=  c_{\gamma\upsilon}  \ee^{\gamma \upsilon t\over 2}}^{\min(p,n-1)   }   c\ee^{-(\ln k)^2}  +  \left[  c(1+ ((2 \times \frac{2}{\gamma \upsilon}+1) \ln k -t )^{\eta}  ) \ee^{-  ( (2 \times \frac{2}{\gamma \upsilon}+1) \ln k -t )}    \right]     
  \\
  &\leq  \sum_{p= c_{\gamma\upsilon}  \ee^{\gamma \upsilon t\over 2} }^{+\infty}   \ee^{t-x} p^{r(\epsilon_0)} \P\left( \min_{j\leq p}(S_j-r(\epsilon_0)\ln j)\geq -x,\, S_p\leq r(\epsilon_0)\ln p+t-x\right)    \times \left[ c \ee^{ -(\frac{\gamma \upsilon}{2}t)^2}  +  c\ee^{-  ( (2 \times \frac{2}{\gamma \upsilon}+\frac{1}{2}) \frac{\gamma\upsilon}{2} t -t )}  \right]
  \end{align*}
  where we have used in the last line the many-to-one Lemma and the trivial inequality  $x\ee^{-(\theta+1)x} \leq c\ee^{-(\theta+\frac{1}{2}) x}$ $\forall x\geq 0$ for some $c>0$. Now by using the inequality $ \P\left( \min_{j\leq p}(S_j-r(\epsilon_0)\ln j)\geq ...\right) \leq \P\left( \min_{j\leq p}(S_j-r(\epsilon_0)\ln j)\geq ... \right)^{1-\frac{3}{2}\frac{\epsilon_0}{2}} $ then  (\ref{estimAid}), we get
  \begin{align*}
 & \P\left( \max_{|z|=n} \sum_{k=   c_{\gamma\upsilon}  \ee^{\gamma \upsilon t\over 2} }^{n-1}   1_{\{  V(z_{|k})\geq \kappa\ln k-x,\,  \exists u\geq z_{|k}^{(l)},\, V(u)<r(\epsilon_0)\ln |u|-x+t+1  \}}  \geq   c_{\gamma\upsilon} \ee^{\gamma \upsilon t\over 2}    ,\, A(x,\epsilon_0) \right)
  \\
  &\leq \sum_{p= c_{\gamma\upsilon}  \ee^{\gamma \upsilon t\over 2} }^{+\infty}   \ee^{t-x} \frac{p^{r(\epsilon_0)}}{p^{\frac{3}{2}-\frac{\epsilon_0}{2}}}\times (1+x)^{\eta} (1+r(\epsilon)\ln p+t)^\eta    \times \left[c \ee^{  -(\frac{\gamma \upsilon}{2}t)^2}  +  c\ee^{-   (t+\frac{1}{2} \frac{\gamma\upsilon}{2} t  )}  \right]
  \end{align*}
 with $\eta= 1-\frac{2}{3}\frac{\epsilon_0}{2}$. As $\frac{3}{2}-\frac{\epsilon_0}{2}-r(\epsilon_0)>1$ we deduce that:
% 
% \begin{eqnarray*}
% &&\leq  \sum_{p=\frac{c_{\gamma\eta}}{2} \ee^{\gamma \eta t\over 2} }^{+\infty} \ee^{-x+t}\frac{p^{r(\epsilon_0)}}{p^\frac{3}{2}}   \sum_{j= \frac{c_{\gamma\eta}}{2} \ee^{\gamma \eta t\over 2}}^{\min(p,n-1)} \left[   (1+(\kappa \ln j-t-r(\epsilon_0)\ln j)^\eta)\ee^{t+r(\epsilon_0)\ln j-(\kappa-1) \ln j} \right.
% \\
% &&\qquad\qquad\qquad\qquad \qquad\qquad \left. + \Q\left( V(B(w_{k+1}))-V(w_k)\geq \ln j\right)  \right],
%\end{eqnarray*}
%where $B(w_{k+1})$ is the child of $w_k$ which is not $w_{k+1}$. The sequence $(V(B(w_{k+1})) -V(w_k))_{k\geq 0}$ is i.i.d, and $\Q\left( V(B(w_{k+1}))-V(w_k)\geq \ln j\right) \leq  c\ee^{-\frac{(\ln j)^2}{2}}$.

\begin{align*}
&\P\left( \max_{|z|=n} \sum_{k=   c_{\gamma\upsilon}  \ee^{\gamma \upsilon t\over 2} }^{n-1}   1_{\{  V(z_{|k})\geq \kappa\ln k-x,\,  \exists u\geq z_{|k}^{(l)},\, V(u)<r(\epsilon_0)\ln |u|-x+t+1  \}}  \geq   c_{\gamma\upsilon} \ee^{\gamma \upsilon t\over 2}    ,\, A(x,\epsilon_0) \right)
 \\
 &\leq (1+x)^\eta\left[  \ee^{-x} t^\eta \ee^{t- t - \frac{\gamma\eta}{4} t} + \ee^{-x} \ee^{t} \ee^{-(\frac{\gamma \eta}{2}t)^2}\right] \leq c (1+x)^{\eta} \ee^{-x} \ee^{-ct},
\end{align*}
which may to conclude because of $(\kappa-2-r(\epsilon_0))\frac{\gamma n}{2}>2$. For some $c>0$.
\hfill$\Box$
\\
\\

\subsection{Proof that $M^{\gamma, \beta}$ is not of bounded variation}\label{infinitevar}
%%%%%%%%%%%%%%%%%%%%%%%%%%%%%%%%%%%%%%%%%%
Here, we show that $M^{\gamma, \beta}$ has infinite variation. Let us suppose that $M^{\gamma, \beta}$ has finite total variation with positive probability. A $0-1$ argument tells us that this occurs in fact with probability $1$ (the reader can check that this property is measurable with respect to the asymptotic sigma algebra generated by the weights $(\Theta_u)_{u\in \bigcup_n \T_{\leq n}}$). The 
 associated total variation function $t \mapsto V^{\gamma, \beta}[0,t]$ is a (random) non decreasing function, which therefore defines a positive measure on $[0,1]$. This measure satisfies the cascading rule
\begin{equation*}
V^{\gamma, \beta}[0,1]= \ee^{-\gamma V(0)}  V^{1,\gamma, \beta}[0,1]+\ee^{-\gamma V(1)}  V^{2,\gamma, \beta}[0,1]
\end{equation*}
where $V^{1,\gamma, \beta}[0,1]$ and $V^{2,\gamma, \beta}[0,1]$ are independent copies of $V^{\gamma, \beta}[0,1]$, independent of $V(0),V(1) $.  This type of distributional equation for the total mass of a random measure has been investigated in \cite[proof of Theorem 1]{Rnew3}, based on a fixed point equation studied in \cite{Liu}. Since the mapping $t\mapsto \ln \E[\sum_{|u|=1}e^{-t\gamma V(u)}]$ vanishes for $t=\frac{1}{\gamma}$ as well as  the derivative at this point,  by \cite{Liu} and the measure extension in \cite{Rnew3}, the mapping $t\mapsto V^{\gamma, \beta}[0,t]$ is   a time changed stable process with parameter $\alpha=\frac{1}{\gamma}$  and the  time change is continuous: continuity of the time change has been proved in \cite{BKNSW} as it is nothing but  the derivative martingale associated to the real part of the branching random walk (we do not write its explicit form as we only need its continuity). This shows that $t\mapsto V^{\gamma, \beta}[0,t]$ cannot be continuous. This is a contradiction with the fact that $M^{\gamma, \beta}$ is continuous, as the total variation of a continuous function is necessarily  continuous.

\section{$4-$ moment of $|M_n^{\gamma,\beta}|$}  \label{secmoments}
%%%%%%%%%%%%%%%%%%%%%%%%%%%%%%%%%%%%%%%%%%%%%%%%%%%%%%%%%%%%%%%%%%%%%%%%%%%%%%%%%%%%

We first give a lemma on the tail behaviour of $|M^{\gamma,\beta}_n|$ which is a direct consequence of the moment estimate given by lemma \ref{momentp2} below:  

\begin{Lemma}
\label{maxondyadique}
Let $\alpha>0$. For any $a,\, x\geq 0,\, y\in \r, k,n\in \N$,
\begin{equation}
\P\left( \max_{|z|=k} \ee^{-\gamma V(z)} 1_{\{ \min_{j\leq k}{V}(z_{|j})\geq -a\}} |M_{k+n}^{\gamma,\beta}(z)| 1_{\{   \forall u\geq z,\, V(u)\geq -x+ \alpha \ln |u|  \}} \geq \ee^{\gamma y}  \right)  \leq  c k^{\alpha(1-4 \gamma)} \ee^{-x+4\gamma(x-y)}  \P_a(\min_{j\leq k}S_j\geq 0)
\end{equation}
\end{Lemma}
\noindent{\it Proof of Lemma \ref{maxondyadique}.} First, as each term is positive, we raise the inequality inside the probability to the power four. Then we observe that $\max_{|z|=k}... \leq \sum_{|z|=k}...$, finally by using the Markov inequality we get:
\begin{eqnarray*}
&&\P\left( \max_{|z|=k} \ee^{-\gamma V(z)}  1_{\{  \min_{j\leq k}{V}(z_{|j})\geq -a\}}  |M_{k+n}^{\gamma,\beta}(z)| 1_{\{   \forall u\geq z,\, V(u)\geq -x+ \alpha \ln |u|  \}} \geq \ee^{\gamma y} \right)
\\
&&\leq \E\left(  \sum_{|z|=k} 1_{\{   \min_{j\leq k}{V}(z_{|j})\geq -a\}} \ee^{-4\gamma V(z)} |M^{\gamma,\beta}_{k+n}(z)|^{4}  1_{\{ \forall u\geq z,\, V(u)-V(z)\geq -x-V(z)+ \alpha \ln |u|  \}}  \right) \ee^{-4\gamma y}
\end{eqnarray*}
Via the branching property and Lemma \ref{momentp2} below, we have
\begin{eqnarray*}
&&\P\left( \max_{|z|=k} \ee^{-\gamma V(z)}  1_{\{  \min_{j\leq k}{V}(z_{|j})\geq -a\}} |M_{k+n}^{\gamma,\beta}(z)| 1_{\{   \forall u\geq z,\, V(u)\geq -x+ \alpha \ln |u|  \}} \geq \ee^{\gamma y} \right)
\\
&&\leq c(p) \E\left(  \sum_{|z|=k}  1_{\{   \min_{j\leq k}{V}(z_{|j})\geq -a\}} \ee^{-4\gamma V(z)} \ee^{4\gamma (x+ V(z)-\alpha \ln k)} \ee^{-(x+V(z)-\alpha \ln k)}   \right) \ee^{-4\gamma y}
\\
&&\leq  c(p) \E\left(  \sum_{|z|=k}  1_{\{  \min_{j\leq k}{V}(z_{|j})\geq -a\}} \ee^{-V(z)} k^{\alpha-4\gamma\alpha}    \right) \ee^{4\gamma (x-y)}\ee^{-x}
\\
&&\leq c(p) k^{\alpha(1-4 \gamma)} \ee^{4\gamma(x-y)} \ee^{-x}  \P_a(   \min_{j\leq k}S_j\geq 0).
\end{eqnarray*}
This completes the proof of Lemma \ref{maxondyadique}.\qed
\medskip

We state now the main moment estimate of the paper:

\begin{Lemma}
	\label{momentp2}
There exists  $ c>0$ such that for any $x\geq 0$, $n\in \N$,
	\begin{equation}
	\E\left(  |M^{\gamma,\beta}_n|^{4} 1_{\{\forall u\in \mathbb{T},\, V(u)\geq -x    \}}   \right)\leq c \ee^{4\gamma x} \ee^{-x} .
	\end{equation}
\end{Lemma}

\noindent{\it Proof of Lemma \ref{momentp2}.} On the set $\{  \forall u\in \mathbb{T},\, V(u)\geq -x\}$, we will bound the expectation of 
\begin{eqnarray} 
\nonumber \label{sum} |M^{\gamma,\beta}_n|^{4}&=& (M^{\gamma,\beta})^2(\overline{M^{\gamma,\beta}})^2
\\
&= &  \sum_{|z_1|=|z_2|=|z_3|=| z_4|=n} \ee^{-\gamma [V(z_1)+V(z_2)+V( z_3)+V( z_4)]} \ee^{i\beta\sqrt{2\ln 2}[X(  z_1)-X(  z_2)+X(  z_3)-X( z_4)]}.
\end{eqnarray}
To carry out this computation, we must look at the different configurations for the genealogical structure of $z_1,z_2,z_3,z_4$. So for any $i,j\in \{1,2,3,4\}$ with $i< j$, we introduce $T(i,j):=\sup{\{ k\in [0,n],\, {z_i}_{|k}={z_j}_{|k}  \}}$, i.e the vertex ${z_i}_{|T(i,j)}={z_j}_{|T(i,j)}$ is the youngest common ancestor of the vertices $z_{i}$ and $z_j$. We also define $T:= \inf_{i,j}T(i,j)$ the generation of the youngest common ancestor of $z_1$, $z_2$, $z_3$ and $z_4$.

Let us decompose the sum (\ref{sum}), according to the different configurations of the $T(i,j)$, i.e let us write
\begin{eqnarray}
\nonumber &&|M^{\gamma,\beta}_n|^{4}= \sum_{  0\leq  s_{ij}\leq n,\, i<j\in \{1,...,4\}} \sum_{|z_1|=...=| z_4|=n  } 1_{\{ \forall i<j\in \{1,...,4\},\, s_{ij}=T(i,j) \}}\times 
\\
\label{sum2}&&   \qquad\qquad\qquad\ee^{-\gamma [V(z_1)+V(z_2)+V( z_3)+V( z_4)]} \ee^{i\beta\sqrt{2\ln 2}[X(  z_1)-X(  z_2)+X(  z_3)-X( z_4)]}
\end{eqnarray}
Write with full details how to bound the expectation of every term of this big sum is a very cumbersome work. However notice that the terms which have the same genealogical structure can be grouped in bundles. Moreover as many of these bundles are very similar or even, because of the symmetry, identical or identical in law, it is easy to convince oneself that it suffices to only treat the following $3$ situations:
\begin{itemize}
	\item 1) $T=n$ (which means that $z_1=z_2=z_3=z_4$). 
	
	\item 2) $T=T(1,2)<T(2,3)<T(3,4)<n$.
	
	\item 3) $T=T(1,3)=T(1,4)$,  $T<T(1,2)<n$ and $T<T(3,4)<n$.
\end{itemize} 

\paragraph{Situation 1)} This is the   simplest one. The restriction of the sum in $\eqref{sum2}$ to the terms for which $T=n$ is simply $S_1:=\sum_{|z|=n}\ee^{-4 \gamma V(z)} $. Note that the imaginary part cancels with itself so it remains only positive terms. On the set $\{\forall u\in \mathbb{T},\, V(u)\geq -x\}$ we thus have
\begin{eqnarray*}
\E(S_1  1_{\{ \forall u\in \mathbb{T},\, V(u)\geq -x \}})&\leq & \E\left(  \sum_{|z|=n}\ee^{-4 \gamma V(z)}1_{\{ \min_{i\leq n} V(u_{|i})\geq -x \}}  \right)
\\
&= & \E\left( \ee^{(1-4\gamma)S_n}1_{\{ \min_{i\leq n} S_i\geq -x \}}\right)
\\
&\leq &  \ee^{4\gamma x}\ee^{-x},
\end{eqnarray*}
which is the desired bound. 

\paragraph{Situation 2)} This case requires to be a bit more careful. Recall that for $r,n\in \N$, $n>r$, $z\in \mathbb{T}$ with $|z|=n$, the vertex ${z}_{|r}$ denotes the ancestor of $z$ in the generation $r$. We now define ${z}_{|r*}$ the child of ${z}_{|r}$ which is not an ancestor of $z$. The restriction of the sum in (\ref{sum2}) to the terms for which $T=T(1,2)<T(2,3)<T(3,4)<n$ can be expressed as 
\begin{eqnarray*}
&&S_2:= \sum_{|z|=n}    \sum_{t=0}^{n-3} \left(  \left\{ \ee^{-\gamma V({z}_{|t*}) +i\beta \sqrt{2\ln 2} X({z}_{|t*})} \times   M_n^{\gamma,\beta}({z}_{|t*})   \right\} \times \sum_{s=t+1}^{n-2}   \left( \left\{ \ee^{-\gamma V({z}_{|s*}) -i\beta \sqrt{2\ln 2} X({z}_{|s*})} \times  \overline{M_n^{\gamma,\beta}({z}_{|s*})  } \right\}  \right. \right.
\\
&&\qquad\qquad\qquad\qquad \qquad \left.  \times    \sum_{r=s+1}^{n-1}   \left\{ \ee^{-\gamma V({z}_{|r*}) +i\beta \sqrt{2\ln 2} X({z}_{|r*})} \times   M_n^{\gamma,\beta}({z}_{|r*})   \right\}\ee^{-\gamma V(z) -i\beta\sqrt{2\ln 2} X(z)} \Big)  \right) .
\end{eqnarray*}
It corresponds to a decomposition of the sum (\ref{sum2}) through the position in the binary tree of the particle denoted by $z_4$ in (\ref{sum2}). Recall that we need to bound the expectation of $S_2$ on the set $\{\forall u\in \mathbb{T},\, V(u)\geq -x\}$. To treat the complex part we will first condition on $\sigma(V(z),z \in \mathbb{T})$, the sigma field generated by the real part. 

To compute this conditional expectation, first observe that for each $|z|=n$, $0\leq t<s<r\leq n-1$, conditionally to $\sigma(V(z),z \in \mathbb{T})$, the  random variables $M_n^{\gamma,\beta}(z_{|r*})$, $\overline{ M_n^{\gamma,\beta}(z_{|s*})}$ and $M_n^{\gamma,\beta}(z_{|t*})$, are mutually independent and independent of others random variables. Moreover their conditional expectations are respectively equal to $2^{-\beta^2(n-t-1)}M_n^{\gamma}(z_{|t*})$, $2^{-\beta^2(n-s-1)}  M_n^{\gamma}(z_{|s*})$ and $2^{-\beta^2(n-r-1)}  M_n^{\gamma}(z_{|r*})$. 

Furthermore note that 
$$\ee^{i\beta \sqrt{2\ln 2} X(z_{|r*})}\ee^{-i\beta \sqrt{2\ln 2}X(z)}= \ee^{i\beta \sqrt{2\ln 2} (X(z_{|r*})-X(z_{|r}))} \ee^{-i\beta \sqrt{2\ln 2}(  X(z)-X(z_{|r} ))}$$ 
is the product of two independent terms whose expectation is equal to $2^{-\beta^2}  2^{-\beta^2 (n-r)}$. Finally by using the same decomposition for $ \ee^{i\beta \sqrt{2\ln 2} X(z_{|t*})} \times \ee^{-i\beta \sqrt{2\ln 2} X(z_{|s*})}$ we get that
\begin{eqnarray*}
&&  \E\left( 1_{\{\forall u\in \mathbb{T},\, V(u)\geq -x  \}} S_2\right)
\\
&&= \E\left(   1_{\{\forall u\in \mathbb{T},\, V(u)\geq -x  \}}  \sum_{|z|=n}  \times \sum_{t=0}^{n-3}  \left\{ 2^{-\beta^2} \ee^{-\gamma V({z}_{|t*})}   2^{-\beta^2(n-t-1)} M_n^{\gamma}({z}_{|t*})   \right\}   \times  \sum_{s=t+1}^{n-2}    \left\{ 2^{-\beta^2(s-(t+1))}\ee^{-\gamma V({z}_{|s*}) }  \right.\right.
\\
&& \qquad\qquad\qquad   \left.  2^{-\beta^2(n-s-1)}  M_n^{\gamma}({z}_{|s*})   \right\}\left. \times  \sum_{r=s+1}^{n-1}   \left\{ 2^{-\beta^2} \ee^{-\gamma V({z}_{|r*}) }  2^{-\beta^2(n-r-1)}  M_n^{\gamma}({z}_{|r*})   \right\}  \ee^{-\gamma V(z)} 2^{-\beta^2 (n-r)}    \right).
\end{eqnarray*}
By reordering the sum, this is equal to
\begin{eqnarray*}
&&=   \sum_{t=0}^{n-3}  \sum_{s=t+1}^{n-2}     \sum_{r=s+1}^{n-1} \sum_{|z|=r} \E\left(   1_{\{\forall u\in \mathbb{T},\, V(u)\geq -x  \}}   \left\{ 2^{-\beta^2} \ee^{-\gamma V({z}_{|t*})}   2^{-\beta^2(n-t-1)} M_n^{\gamma}({z}_{|t*})   \right\}   \times   \left\{ 2^{-\beta^2(t+1-s)}\ee^{-\gamma V({z}_{|s*}) }  \right.\right.
\\
&&   \left.   2^{-\beta^2(n-s-1)}  M_n^{\gamma}({z}_{|s*})   \right\}\left. \times     \sum_{|z'|=n,\, z'\geq z}  2^{-\beta^2} \ee^{-\gamma V({z'}_{|r*}) } M_n^{\gamma}({z'}_{|r*})   2^{-\beta^2(n-r-1)}    \ee^{-\gamma V(z')} 2^{-\beta^2 (n-r)}    \right).
\end{eqnarray*}
We are left with positive terms only. So by estimating $1_{\{ \forall u\in \mathbb{T},\, V(u)\geq -x  \}} $ with $  1_{\{ \min_{i\leq |z|}V(z_{|i})\geq -x \}}$ for any $z\in \mathbb{T}$ such that $|z|=r$, we can use the branching property for the real part of the branching random walk. Indeed for any $z\in \mathbb{T}$, let $\mathcal{F}^c(z)=\sigma ( V(u),\, \text{such that } u \text{ is not a descendent of } z  )$, by the branching property we have 
$$2^{ -\beta^{2}(n-s-1) }\E\left( M_n^{\gamma}(z_{|s*})\big| \mathcal{F}^c(z_{|s*}) \right)=2^{((1-\gamma)^2-\beta^2)(n-s-1)}=1,$$
where we observe that $\gamma+\beta=1$ implies $(1-\gamma)^2-\beta^2=0$. Similarly it is also plain to check that for any $|z|=r$, 
\begin{equation*}
\E\left(   \sum_{|z'|=n,\, z'\geq z}  2^{-\beta^2} \ee^{-\gamma V({z'}_{|r*}) } M_n^{\gamma}({z'}_{|r*})   2^{-\beta^2(n-r-1)}    \ee^{-\gamma V(z')} 2^{-\beta^2 (n-r)} |\mathcal{F}_r \right)=2^{-1} \ee^{-\gamma V(z)},
\end{equation*}
with $\mathcal{F}_r:= \sigma(V(u),\, |u|\leq r)$. Combining these arguments we get
\begin{eqnarray*}
  \E\left( 1_{\{\forall u\in \mathbb{T},\, V(u)\geq -x  \}} S_2\right)  
  &\leq&  c\sum_{t=0}^{n-3} \sum_{s=t+1}^{n-2} \sum_{r=s+1}^{n-1}  \E\left(   \sum_{|z|=r}       \ee^{-\gamma V({z}_{|t})}           \ee^{-\gamma V({z}_{|s}) }      \ee^{-2\gamma V({z}) }     1_{\{   \min_{j\leq r}V(z_{|j}) \geq -x\}}    \right) 
  \\
  &= & c \sum_{t=0}^{n-3} \sum_{s=t+1}^{n-2} \sum_{r=s+1}^{n-1}  \E\left(       \ee^{-\gamma S_t}           \ee^{-\gamma S_s }      \ee^{(1-2\gamma) S_r }     1_{\{   \min_{j\leq r}S_j \geq -x\}}    \right) ,
\end{eqnarray*}
where the equality stems from the the many to one Lemma (\ref{manytoone}). Finally by using three times the inequality (\ref{estimAid2}) (with first $\kappa:= 2\gamma-1$, then twice $\kappa=\gamma$) we get
\begin{eqnarray*}
&& \E\left( 1_{\{\forall u\in \mathbb{T},\, V(u)\geq -x  \}} S_2\right)
\\
&&\leq c   \sum_{t=0}^{n-3} \sum_{s=t+1}^{n-2}    \E\left(       \ee^{-\gamma S_t}      \times     \ee^{(1-3\gamma) S_s }  1_{\{   \min_{j\leq s}S_j \geq -x\}}        \sum_{r=s+1}^{n-1}\E\left( \ee^{(1-2\gamma)S_{r-s}}   1_{\{   \min_{j\leq r-s}S_j \geq -x-S_s\}}       \right)     \right)
\\
&&\leq c' \sum_{t=0}^{n-3}\sum_{s=t+1}^{n-2}\E\left(  \ee^{-\gamma S_t} \ee^{-\gamma S_s} 1_{\{ \min_{j\leq s}S_j\geq -x  \}}  \right) \ee^{(2\gamma-1)x}
  \\
  &&\leq   c''\sum_{t=0}^{n-3} \E\left(       \ee^{-\gamma S_t}             1_{\{   \min_{j\leq t}S_j \geq -x\}} \right) \ee^{-x+3\gamma x}\leq c''\ee^{-x+4\gamma x}.
\end{eqnarray*}
We are done with the study of the situation 2).

\paragraph{Situation 3)} This is the trickiest case. The restriction of the sum in $\eqref{sum2}$ to the terms for which $T=T(1,3)=T(1,4)$,  $T<T(1,2)<n$ and $T<T(3,4)<n$ can be rewritten 
\begin{eqnarray}
	\label{thesum}
&&S_3:= \sum_{t=0}^{n-2}\sum_{s=t+1}^{n-1} \sum_{s'=t+1}^{n-1}\sum_{|z_1|=...=| z_4|=n  } 1_{\{ T(1,3)=T(1,4)=t,  T(1,2)=s,\, T(3,4)=s' \}}\times 
\\
\nonumber &&   \qquad\qquad\qquad\ee^{-\gamma [V(z_1)+V(z_2)+V( z_3)+V( z_4)]} \ee^{i\beta\sqrt{2\ln 2}[X(  z_1)-X(  z_2)+X(  z_3)-X( z_4)]}.
\end{eqnarray}
We denote by $a$ the youngest ancestor of $z_1,z_2,z_3$ and $z_4$, $b$ the youngest ancestor of $z_1$ and $z_2$ and $c$ this one of $z_3$ and $z_4$. To express precisely the sum \eqref{thesum} we need to specify whether $z_1$ and $z_2$ are descendant of $a^{(l)}$ or $a^{(r)}$ then whether $z_1$ is a descendant of $b^{(l)}$ or $b^{(r)}$ and whether $z_3$ is descendant of $c^{(l)}$ or $c^{(r)}$. It gives $8$ different configurations $e\in \{1,...,8\}$ and leads to consider the random variable $S^{(e)}_3$ defined as the sum of the terms for which $z_1$, $z_2$, $z_3$ and $z_4$ are in the configuration $e$. It is easy to see that these $(S^{(e)}_3)_{e\in \{1,...,8\}}$ have the same law. Then we just have to   compute $S^{(1)}_3$, the sum of the terms for which $z_1$ and $z_2$ are descendant of $a^{(l)}$, $z_1$ is a descendant of $b^{(l)}$ and $z_3$ a descendant of $c^{(l)}$.

% By symmetry we see that the sum of the terms for which $z_1$ and $z_2$ are descendant of $a^{(l)}$ (the left child of $a$)  is equal to the sum of the term for which $z_1$ and $z_2$ are descendant of $a^{(r)}$. 
%
%Similarly the sum of the terms for which $z_1$ is a descendant of $b^{(l)}$ and $z_3$ is a descendant of $c^{(l)}$ is the complex conjugate of the sum of the terms for which $z_1$ is a descendant of $b^{(r)}$ and $z_3$ is a descendant of $c^{(r)}$. And of course it remains the case where $z_1$ is a descendant of $b^{(l)}$ 
%
Finally the expectation on the set $\{\forall u\in \mathbb{T},\, V(u)\geq -x\}$ of the complicated sum in (\ref{sum2}), restricted to the terms for which $T=T(1,3)=T(1,4)$,  $T<T(1,2)<n$ and $T<T(3,4)<n$, is equal to 
\begin{eqnarray*}
&& \E( 1_{\{  \forall u\in \mathbb{T},\, V(u)\geq -x \}}S_3 ) =8\E\Big[ 1_{\{  \forall u\in \mathbb{T},\, V(u)\geq -x \}} \sum_{t=0}^{n-2}\sum_{|a|=t}\ee^{-4\gamma V(a)}  \times 
\\
&& \Big\{ \sum_{s=t+1}^{n-1} \sum_{|b|=s,\, b\geq a^{(l)}} \ee^{-2\gamma (V(b)-V(a) )}\times \left( \ee^{-\gamma (V(b^{(l)})-V(b))+i\beta \sqrt{2\ln 2}(X(b^{(l)})-X(b))} M_n^{\gamma,\beta}(b^{(l)})  \right) \times 
\\
&&\qquad\qquad\qquad  \left(  \ee^{-\gamma (V(b^{(r)})-V(b))-i\beta \sqrt{2\ln 2}(X(b^{(r)})-X(b))} \overline{M_n^{\gamma,\beta}(b^{(r)})  } \right)  {\Big\} }\times
\\
&& {\Big\{ }\sum_{s'=t+1}^{n-1} \sum_{|c|=s',\, c\geq a^{(r)}} \ee^{-2\gamma (V(c)-V(a) )}\times \left( \ee^{-\gamma (V(c^{(l)})-V(c))+i\beta \sqrt{2\ln 2}(X(c^{(l)})-X(c))} M_n^{\gamma,\beta}(c^{(l)})  \right) \times 
\\
&&\qquad\qquad\qquad  \left(  \ee^{-\gamma (V(c^{(r)})-V(c))-i\beta \sqrt{2\ln 2}(X(c^{(r)})-X(c))} \overline{M_n^{\gamma,\beta}(c^{(r)})  } \right)  {\Big\} }  \Big].
\end{eqnarray*}
Now by taking the conditional expectation with respect to $\sigma(V(z),z\in \mathbb{T})$, the sigma-field generated by the real part, and by using the branching property (for the complex part as in situation 2)) we get
\begin{eqnarray*}
	&& \E( 1_{\{  \forall u\in \mathbb{T},\, V(u)\geq -x \}}S_3 ) =8\E\Big[ 1_{\{  \forall u\in \mathbb{T},\, V(u)\geq -x \}} \sum_{t=0}^{n-2}\sum_{|a|=t}\ee^{-4\gamma V(a)}  \times 
	\\
	&& \Big\{ \sum_{s=t+1}^{n-1} \sum_{|b|=s,\, b\geq a^{(l)}} \ee^{-2\gamma (V(b)-V(a) )}\times \left( 2^{-\beta^{2}} \ee^{-\gamma (V(b^{(l)})-V(b))} 2^{-\beta^2 (n-s-1)}M_n^{\gamma}(b^{(l)})  \right) \times 
	\\
	&& \left( 2^{-\beta^{2}} \ee^{-\gamma (V(b^{(r)})-V(b))} 2^{-\beta^{2}(n-s-1)}M_n^{\gamma}(b^{(r)})   \right)  {\Big\} }\times
	\\
	&& {\Big\{ }\sum_{s'=t+1}^{n-1} \sum_{|c|=s',\, c\geq a^{(r)}} \ee^{-2\gamma (V(c)-V(a) )}\times \left( 2^{-\beta^2}\ee^{-\gamma (V(c^{(l)})-V(c))} 2^{-\beta^{2}(n-s'-1)} M_n^{\gamma}(c^{(l)})  \right) \times 
	\\
	&& \left( 2^{-\beta^2} \ee^{-\gamma (V(c^{(r)})-V(c))} 2^{-\beta^{2}(n-s'-1)}  M_n^{\gamma}(c^{(r)})   \right)  {\Big\} }  \Big].
\end{eqnarray*}
Again, we are left with positive terms only.  So using the estimates $1_{\{ \forall u\in \mathbb{T},\, V(u)\geq -x  \}}  \leq 1_{\{ \min_{i\leq |b|}V(b_{|i})\geq -x \}}$ and $ 1_{\{ \forall u\in \mathbb{T},\, V(u)\geq -x  \}}  \leq 1_{\{ \min_{i\leq |c|}V(c_{|i})\geq -x \}} $ for any $b,c\in \mathbb{T}$ with $|b|=s$, $|c|=s'$ and the branching property, we get
%\begin{eqnarray*}
%	&& \E( 1_{\{  \forall u\in \mathbb{T},\, V(u)\geq -x \}}S_2 ) =8\E\Big[    \sum_{t=1}^{n-2}\sum_{|a|=t}\ee^{-4\gamma V(a)}  \times 
%	\\
%	&& \Big\{ \sum_{s=t+1}^{n-1} \sum_{|b|=s,\, b\geq a^{(l)}} \ee^{-2\gamma (V(b)-V(a) )}   1_{\{ \min_{j\leq s} V(b_j) \geq -x \}}  \times \left( 2^{-\beta^{2}} \ee^{-\gamma (V(b^{(l)})-V(b))} 2^{-\beta^2 (n-s-1)}M^{\gamma}(b^{(l)})  \right) \times 
%	\\
%	&& \left( 2^{-\beta^{2}} \ee^{-\gamma (V(b^{(r)})-V(b))} 2^{-\beta^{2}(n-s-1)}M^{\gamma}(b^{(r)})   \right)  {\Big\} }\times
%	\\
%	&& {\Big\{ }\sum_{s'=t+1}^{n-1} \sum_{|c|=s',\, c\geq a^{(r)}} \ee^{-2\gamma (V(c)-V(a) )}  1_{\{ \min_{j\leq s} V(c_j) \geq -x \}}   \times \left( 2^{-\beta^2}\ee^{-\gamma (V(c^{(l)})-V(c))} 2^{-\beta^{2}(n-s'-1)} M^{\gamma}(b^{(l)})  \right) \times 
%	\\
%	&& \left( 2^{-\beta^2} \ee^{-\gamma (V(c^{(r)})-V(c))} 2^{-\beta^{2}(n-s'-1)}  M^{\gamma}(c^{(r)})   \right)  {\Big\} }  \Big].
%\end{eqnarray*}
%Now by applying the branching property, we get
\begin{eqnarray*}
%	&& \E( 1_{\{  \forall u\in \mathbb{T},\, V(u)\geq -x \}}S_3 )
%	\\
%	&& =8\E\Big[    \sum_{t=0}^{n-2}\sum_{|a|=t}\ee^{-4\gamma V(a)}  \times 
%	\\
%	&& \Big\{ \sum_{s=t+1}^{n-1} \sum_{|b|=s,\, b\geq a^{(l)}} \ee^{-2\gamma (V(b)-V(a) )}   1_{\{ \min_{j\leq s} V(b_{|j}) \geq -x \}}  \times \left( 2^{-\beta^{2}} 2^{\gamma^2-2\gamma} 2^{-\beta^2 (n-s-1)}  2^{(1+\gamma^2-2\gamma)(n-s-1)} \right) \times 
%	\\
%	&& \left( 2^{-\beta^{2}}  2^{\gamma^2-2\gamma} 2^{-\beta^{2}(n-s-1)} 2^{(1+\gamma^2-2\gamma)(n-s-1)}  \right)  {\Big\} }\times
%	\\
%	&& {\Big\{ }\sum_{s'=t+1}^{n-1} \sum_{|c|=s',\, c\geq a^{(r)}} \ee^{-2\gamma (V(c)-V(a) )}  1_{\{ \min_{j\leq s} V(c_{|j}) \geq -x \}}   \times \left( 2^{-\beta^2} 2^{\gamma^2-2\gamma} 2^{-\beta^{2}(n-s'-1)}  2^{(1+\gamma^2-2\gamma)(n-s'-1)}  \right) \times 
%	\\
%	&& \left( 2^{-\beta^2}  2^{\gamma^2-2\gamma} 2^{-\beta^{2}(n-s'-1)}   2^{(1+\gamma^2-2\gamma)(n-s'-1)}  \right)  {\Big\} }  \Big].
%\\
	&& \E( 1_{\{  \forall u\in \mathbb{T},\, V(u)\geq -x \}}S_3 )	\leq \frac{1}{2}\E\Big[    \sum_{t=0}^{n-2}\sum_{|a|=t}\ee^{-4\gamma V(a)} \times \Big\{ \sum_{s=t+1}^{n-1} \sum_{|b|=s,\, b\geq a^{(l)}} \ee^{-2\gamma (V(b)-V(a) )}   1_{\{ \min_{j\leq s} V(b_{|j}) \geq -x \}}   \Big\} 
	\\
	&&\qquad\qquad\qquad\qquad\qquad \qquad \qquad  \Big\{ \sum_{s'=t+1}^{n-1} \sum_{|c|=s',\, c\geq a^{(r)}} \ee^{-2\gamma (V(c)-V(a) )}  1_{\{ \min_{j\leq s} V(c_{|j}) \geq -x \}}    \Big\}   \Big],
\end{eqnarray*}
where we have equalities similar to $\E\left( 2^{-\beta^2(n-s'-1)} M^\gamma_n(c^{(r)})\right)= \E\left( 2^{-\beta^2(n-s'-1)} M^\gamma_n(c^{(r)})\right)=1 $. Now by using the many-to-one Lemma (\ref{manytoone}) for the sub-trees $\{ b\in \mathbb{T},\, b\geq a^{(l)}  \}$ and   $\{ c\in \mathbb{T},\, a\geq a^{(r)}  \}$, we get
\begin{eqnarray*}
&& \E( 1_{\{  \forall u\in \mathbb{T},\, V(u)\geq -x \}}S_3 )   \leq  \frac{1}{2}\E\Big[    \sum_{t=0}^{n-2}\sum_{|a|=t}\ee^{-4\gamma V(a)}  1_{\{ \min_{j\leq t} V(a_{|j}) \geq -x \}} \times  
 	\\
 	&&\qquad  \Big\{ \sum_{s=t+1}^{n-1} \E\left( \ee^{(1-2\gamma) S_{s'-t}}1_{\{ \min_{j\leq s'-t}S_j\geq -x-V(a)  \}}  \right)  \Big\} \times \Big\{    \sum_{s'=t+1}^{n-1} \E\left( \ee^{(1-2\gamma) S_{s'-t}}1_{\{ \min_{j\leq s'-t}S_j\geq -x-V(a)  \}}  \right)   \Big\}\Big].
\end{eqnarray*}
Finally by the inequality \eqref{estimAid2} (with $\kappa:= 2\gamma-1$), it follows that
\begin{eqnarray*}
 \E( 1_{\{  \forall u\in \mathbb{T},\, V(u)\geq -x \}}S_3 ) & \leq & c\E\Big[    \sum_{t=0}^{n-2}\sum_{|a|=t}\ee^{-4\gamma V(a)}  1_{\{ \min_{j\leq t} V(a_{|j}) \geq -x \}}    \ee^{2(1-2\gamma) (-x-V(a))} \Big]
\\
&\leq & c   \ee^{-2x +4\gamma x} \E\Big[    \sum_{t=0}^{n-2} \ee^{- S_t}   1_{\{ \min_{j\leq t} S_j \geq -x \}}     \Big]
\\
&\leq & c'\ee^{-x+4\gamma x},
\end{eqnarray*}
where in the last line we have used another time the inequality (\ref{estimAid2}) (with $\kappa=1$). This achieves the study of the situation 3) and thus the proof of Lemma \ref{momentp2}.\hfill$\Box$

\end{document}